\documentclass[12pt]{amsart}
\usepackage{amssymb,amsmath}
\theoremstyle{plain}
\newtheorem{thm}{Theorem}[section]
\newtheorem{cor}[thm]{Corollary}
\newtheorem{prop}[thm]{Proposition}
\newtheorem{lem}[thm]{Lemma}

\theoremstyle{definition}
\newtheorem{rem}[thm]{Remark}
\newtheorem{defn}[thm]{Definition}

\newtheorem{ques}[thm]{Question}

\newcommand{\Prf}{\noindent\textbf{Proof.\ }}
\renewcommand{\qed}{\hfill \vrule height5pt width5pt depth1pt \vspace{+2.00ex}}
\newcommand{\upqed}{\vspace{-2.5\baselineskip}\newline\hbox{}\qed}
\newcommand{\dlim}{\displaystyle\lim\limits}

\newcommand{\ep}{\varepsilon}
\renewcommand{\phi}{\varphi}
\newcommand{\lip}{\langle}
\newcommand{\rip}{\rangle}
\newcommand{\ip}[1]{\lip #1 \rip}

\newcommand{\ol}{\overline}
\newcommand{\one}{\boldsymbol{1}}
\newcommand{\re}{\operatorname{Re}}

\newenvironment{spmatrix}{\left(\begin{smallmatrix}}{\end{smallmatrix}\right)}

\newcommand{\upchi}{{\raise.35ex\hbox{$\chi$}}}
\newcommand{\wot}{\textsc{wot}}
\newcommand{\sotsum}{\textsc{sot--}\!\!\sum}
\newcommand{\wotsum}{\textsc{wot--}\!\!\sum}
\newcommand{\AND}{\text{ and }}
\newcommand{\FOR}{\text{ for }}
\newcommand{\FORAL}{\text{ for all }}

\newcommand{\qand}{\quad\text{and}\quad}

\newcommand{\qforal}{\quad\text{for all}\quad}
\newcommand{\qif}{\quad\text{if}\quad}

\newcommand{\Alg}{\operatorname{Alg}}
\newcommand{\botimes}{\operatorname{\bar{\otimes}}}

\newcommand{\diag}{\operatorname{diag}}
\newcommand{\Dim}{\operatorname{dim}}
\newcommand{\dist}{\operatorname{dist}}
\newcommand{\id}{\operatorname{id}}
\newcommand{\Lat}{\operatorname{Lat}}
\newcommand{\ran}{\operatorname{Ran}}
\newcommand{\rank}{\operatorname{rank}}
\newcommand{\Ref}{\operatorname{Ref}}
\newcommand{\spn}{\operatorname{span}}

\newcommand{\bC}{{\mathbb{C}}}

  \newcommand{\A}{{\mathcal{A}}}
  \newcommand{\B}{{\mathcal{B}}}
  \newcommand{\C}{{\mathcal{C}}}
  \newcommand{\D}{{\mathcal{D}}}

  \newcommand{\G}{{\mathcal{G}}}
\renewcommand{\H}{{\mathcal{H}}}
  \newcommand{\I}{{\mathcal{I}}}
  \newcommand{\J}{{\mathcal{J}}}
  \newcommand{\K}{{\mathcal{K}}}  

  \newcommand{\M}{{\mathcal{M}}}
  \newcommand{\N}{{\mathcal{N}}}

\renewcommand{\P}{{\mathcal{P}}}
  \newcommand{\Q}{{\mathcal{Q}}}
  
\renewcommand{\S}{{\mathcal{S}}}
  \newcommand{\T}{{\mathcal{T}}}
  \newcommand{\U}{{\mathcal{U}}}

\newcommand{\fA}{{\mathfrak{A}}}

\newcommand{\fC}{{\mathfrak{C}}}
\newcommand{\fD}{{\mathfrak{D}}}

\newcommand{\fM}{{\mathfrak{M}}}

\newcommand{\fS}{{\mathfrak{S}}}

\newcommand{\fX}{{\mathfrak{X}}}

\begin{document}

\title{1-Hyperreflexivity and Complete Hyperreflexivity}
%
\author[K.R.Davidson]{Kenneth R. Davidson}
\address{Pure Math.\ Dept.\\U. Waterloo\\Waterloo, ON\;
N2L--3G1\\CANADA}
\email{krdavids@uwaterloo.ca \vspace{-2ex}}
\author[R.H.Levene]{Rupert H. Levene}
\email{rlevene@uwaterloo.ca \vspace{-2ex}}
\date{}
\begin{abstract}
The subspaces and subalgebras of $\B(\H)$ which are hyperreflexive with
constant 1 are completely classified.  It is shown that there are
1-hyper\-reflexive subspaces for which the complete hyperreflexivity
constant is strictly greater than 1. The constants
for $\bC T \otimes \B(\H)$ are analyzed in detail.
\subsection*{Keywords} $1$-hyper\-reflexivity, complete
  hyper\-reflexivity, distance formula
\end{abstract}
\maketitle

The study of invariant subspaces and the notion of reflexivity
plays a central role in operator theory.
The  quantitative notion of hyperreflexivity is a significant
strengthening of reflexivity.
And when this property holds, there are important ramifications.
This is best seen in the theory of nest algebras, where one
obtains a precise distance formula \cite{Arv};
and for a von Neumann algebra, where hyperreflexivity is equivalent to
the vanishing of a certain cohomology group \cite{Chr2}.

Until recently, the collection of known hyperreflexive algebras has been
quite limited.  In addition to nest algebras and most von Neumann algebras
(excluding those with certain intractible type II$_1$ commutants,
where the problem remains open), there were not many others.
The Toeplitz algebra \cite{D2} and certain free semigroup algebras including
the so called noncommutative analytic Toeplitz algebras \cite{DP1,DKP} are
hyperreflexive.
However the first author and S. Power constructed CSL algebras which are not
hyperreflexive \cite{DPo}. 
The class of hyperreflexive algebras was significantly expanded by Bercovici
\cite{Ber} who found general properties which imply hyperreflexivity.
Jaeck and Power \cite{JP} have combined these results to show that the
free semigroupoid algebra associated to any finite directed graph is
hyperreflexive.

These notions make perfect sense for subspaces as well as algebras.
Loginov and Shulman \cite{LS} reformulated reflexivity in this context.
This was done for hyperreflexivity by Larson \cite{Lar}.
See Hadwin \cite{Had1,Had2} for a quite general view of these issues.
A very recent theorem of M\"uller and Ptak \cite{MP} shows that every finite
dimensional reflexive subspace is hyperreflexive, a surprisingly
difficult result.

One focus of this paper is the case where one obtains an exact distance
formula.
We call a subspace 1-hyper\-reflexive if this holds, namely
\[
 \dist(T, \S) = \sup_{\|x\|=1} \dist(Tx,\S x) 
 \qforal T \in \B(\H) .
\]
This can be reformulated as an interchange of sups and infs:
\[
 \inf_{S \in \S} \sup_{\|x\|=1} \| Tx - Sx \| =
 \sup_{\|x\|=1} \inf_{S \in \S} \| Tx - Sx \| .
\]
We will classify these spaces.

The second focus of this paper is the notion of complete hyperreflexivity,
namely the hyperreflexivity of $\S \botimes \B(\H)$, the \wot-closed
spatial tensor product.
It is an open question whether hyperreflexivity of $\S$ implies complete
hyperreflexivity.
In the case of known examples such as nest algebras and von Neumann
algebras, the proofs yields the same constant for the complete case
as for the algebra itself.
The same is true for the Toeplitz algebra, free semigroup algebras and
algebras handled by Bercovici's Theorem.
We will produce examples of hyperreflexive subspaces which are completely
hyperreflexive but for which the constant increases.
Indeed, any one-dimensional subspace $\bC T$ where $\rank T \ge2$ will be
1-hyper\-reflexive but not completely 1-hyper\-reflexive; but it will 
have complete hyperreflexive constant no greater than 4.

\section{Setting the Stage}

Recall that subalgebra $\fA$ of $\B(\H)$ is reflexive if
\[ 
  \fA = \Alg(\Lat \fA) = \{ T \in \B(\H): TP = PTP \FORAL P \in \Lat \fA \} ,
\]
and $\fA$ is hyperreflexive 
if there is a constant $C$ so that for all $T \in \B(\H)$, 
\[ 
 \dist(T,\fA) \le C \sup\{ \|P^\perp T P \| : P \in \Lat \fA \} .
\]
The inequality 
\[
 \sup\{ \|P^\perp T P \| : P \in \Lat \fA \} \le \dist(T,\fA)
 \]
is elementary.
In the same vein, a subspace $\S$ of $\B(\H)$ is reflexive if
\[
 \S = \Ref(\S) = \{ T \in \B(\H): Tx \in \ol{\S x} \FORAL x \in \H \} ,
\]
and $\S$ is hyperreflexive 
if there is a constant $C$ so that for all $T \in \B(\H)$, 
\[ 
 \dist(T,\S) \le C \sup\{ \| P_{\S x}^\perp Tx\| : \|x\|=1 \} .
\]
The optimal constant $\kappa_\S$ is called the \textit{distance constant}.
We say that $\S$ is \textit{$1$-hyper\-re\-flex\-ive} if $\kappa_\S=1$.
We will write 
\[ \beta_\S(T) = \sup\{ \| P_{\S x}^\perp Tx\| : \|x\|=1 \} .\]

One trivial observation is worth recording: \textit{if $\S$ is hyperreflexive
with constant $C$, then so is $U\S V$ where $U$ and $V$ are any
unitary operators.}

\smallskip
One purpose of this paper is to describe all subspaces which are
1-hyper\-re\-flex\-ive.  There are three known classes of algebras
with distance constant 1:

\smallskip\noindent
\textbf{A1.} Nest algebras. Arveson \cite{Arv} (see \cite[Theorem~9.5]{D1}). 

\smallskip\noindent
\textbf{A2.} $\bC I$. Stampfli \cite{St} (see \cite[Theorem~9.15]{D1}).

\smallskip\noindent
\textbf{A3.} $\B(\H_1) \oplus \B(\H_2)$ for any Hilbert spaces
$\H_1$ and $\H_2$. \cite{KL} 
\medskip

Recall that a nest is a chain $\N$ of subspaces of a Hilbert space $\H$
containing both $0$ and $\H$ which is complete with respect to intersections
and closed spans.
The corresponding nest algebra $\T(\N)$ consists of all operators leaving the
nest invariant.
Thus it is reflexive by definition.
The 1-hyperreflexivity of nest algebras is known as the Arveson distance
formula.  It plays a central role in the theory.
We refer the reader to \cite{D1} for more information about nest algebras. 

Stampfli shows that $\dist(T, \bC I) = \frac12 \|\delta_T\|$
where $\delta_T$ is the inner derivation $\delta_T(A) = AT-TA$.
He accomplishes this by proving that if $\dist(T, \bC I) = \|T\| = 1$, then
there is a sequence $x_n$ of unit vectors so that 
\[
 \lim_{n\to\infty} \| Tx_n\|=1 \qand \lim_{n\to\infty} \ip{Tx_n,x_n} = 0 .
\]
Then setting $P_n = x_nx_n^*$ yields 
\[
 \beta_{\bC I}(T) \ge \sup_{n\ge1} \| P_n^\perp TP_n\| = 1 = \|T\| .
\] 

Example A3 cannot be extended to the direct sum of three copies of $\B(\H)$
because even the $3\times3$ diagonal algebra is not 1-hyper\-reflexive.
Indeed it has constant $\sqrt{3/2}$ \cite{DO}.
As we will see below, the proper generalization of this example is
that the space of block off-diagonal operators is 1-hyper\-reflexive.

When we expand our view to subspaces, these examples become:

\smallskip\noindent
\textbf{S1.} Nest bimodules. If $\M$ and $\N$ are nests and
$\theta$ is an order preserving map of $\N$ into $\M$, then the
\wot-closed $\T(\M)$--$\T(\N)$ bimodule 
\[
 \fX(\theta) :=
 \{ T \in \B(\H) : TN \subset \theta(N) \qforal N \in \N \}
\]
 is 1-hyper\-reflexive. 

\smallskip\noindent
\textbf{S2.} $\bC T$ for $T$ an arbitrary operator. Magajna \cite{Mg}.

\smallskip\noindent
\textbf{S3.} Let $\P = \{ P_i : i \in \I \}$ and 
$\Q = \{ Q_i : i \in \I \}$ be partitions of the identity
of $\H_1$ and $\H_2$ respectively.
Then the subspace
\[
 \fX := \{ T \in \B(\H) : Q_iTP_i= 0 \qforal i \in \I \} 
\]
 is 1-hyper\-reflexive.
\medskip

Magajna \cite{Mg} generalizes Stampfli's result in a straighforward
way.  See the remark following Proposition~\ref{P:dist1}.

\smallbreak

If $\fX$ is a $\T(\M)$--$\T(\N)$ nest bimodule which is \wot-closed,
then there is a unique left continuous order preserving map $\theta$ of
$\N$ into $\M$ so that $\fX = \fX(\theta)$ \cite{EP} (see
\cite[Theorem 15.14]{D1}). 
In this paper, all nest bimodules will be \wot-closed;
so we will just call them nest bimodules.
The distance formula for nest bimodules is a routine adaptation
of Power's proof \cite{Po} of Arveson's distance formula.

We wish to isolate part of the ``diagonal'' of a nest~bimodule $\fX$.
It is convenient to describe this by identifying a certain smaller
nest bimodule~$\fX_0$. Consider a finite or countable collection of
elements $\{N_i : i\in \I\}$ in $\N$ such that for every $i$,
$N_i^+\neq N_i$ and $\theta(N_i^+)=M_i^+$ is a successor in~$\M$, and
the restriction of~$\theta$ to this collection is injective. Let
$B_i=N_i^+-N_i$ and $A_i=M_i^+-M_i$ be the corresponding atoms of $\N$
and $\M$ respectively. Define
\[
 \theta_0(N) = \begin{cases} 
                M_i &\qif N = N_i^+,\ i \in \I\\ 
                \theta(N) &\quad\text{otherwise}
               \end{cases} .
\]
This determines a bimodule~$\fX_0$. Moreover, 
\[
\fX=\fX_0+\sum_{i\in\I}\strut^{\!\oplus} A_i \B(\H_1,\H_2)B_i
\]
and we refer to $\sum_{i\in\I}^\oplus A_i \B(\H_1,\H_2)B_i$ as the
diagonal determined by the $\T(\M)$--$\T(\N)$ \textit{bimodule pair}
$(\fX,\fX_0)$. We also write $\Delta(\fX,\fX_0)$ for the set 
$\{(A_i,B_i):i\in \I\}$.

The reason that this definition is convenient is that the notion of an
atom of~$\fX$ is in part determined by the choice of the nests $\M$
and $\N$, and is not intrinsic to~$\fX$. For example, suppose that
$\fX=\B(\H_1,\H_2)$. Then we should choose any proper projections
$A$ and $B$. Consider the nests $\M=\{0,A^\perp\H_2, \H_2\}$ and
$\N=\{0,B\H_1,\H_1\}$. Then $\fX$ is a $\T(\M)$--$\T(\N)$ bimodule
with $\theta(0)=0$ and
$\theta(B\H_1)=\theta(\H_1)=\H_2$. But $A\B(\H_1,\H_2)B$ becomes the
diagonal if we choose $\fX_0$ to be the $\T(\M)$--$\T(\N)$ bimodule
with $\theta(0)=0$, $\theta(B\H_1)=A^\perp\H_2$ and
$\theta(\H_1)=\H_2$.

\medbreak

In the next section, we present two constructions of new 
1-hyper\-re\-flex\-ive subspaces.
Lemma~\ref{tri_const} shows that one can replace atoms of a nest bimodule
with a one-dimensional subspace.
Lemma~\ref{diag_const} shows that in examples of type S3, one can replace
the zero diagonal entries with subspaces formed by Lemma~\ref{tri_const}.
Our goal is to show that every $1$-hyperreflexive subspace is obtained in
this manner.

\pagebreak[3]
\begin{thm}\label{T:1hyp}
Let $\S$ be a \wot-closed subspace of $\B(\H_1,\H_2)$.
Then $\S$ is $1$-hyper\-re\-flex\-ive if and only if there are partitions of
the identity $\C = \{ C_j : j \in \J \}$ and  $\D = \{ D_j : j \in \J \}$
of $\H_2$ and $\H_1$ respectively and for each $j \in \J$, there are
subspaces $\fX_j$ of $C_j \B(\H_1,\H_2) D_j$ obtained from the
construction of Lemma~$\ref{tri_const}$ so that
\[
 \S = \{ T \in \B(\H) : C_jTD_j \in \fX_j \qforal j \in \J \}.
\]
\end{thm}

This will be proven in section~\ref{S:noncomm}.

\section{1-Hyperreflexivity}

In order to complete the list of 1-hyper\-reflexive subspaces, we need
two basic constructions.
The first is more surprising.
It says that atoms of a nest bimodule may be replaced by 
1-dimensional subspaces.

\begin{lem}\label{tri_const}
Let $(\fX,\fX_0)$ be a $\T(\M)$--$\T(\N)$ bimodule pair with diagonal 
$\Delta(\fX,\fX_0)=\{(A_i,B_i) : i \in \I \}$.  
Select operators 
$X_i \in \B(B_i\H_1,A_i\H_2)$\
and define a subspace 
\[
 \S = \{X \in \fX : A_i X B_i \in \bC X_i \FOR i \in \I \} .
\]
Then $\S$ is $1$-hyperreflexive.
\end{lem}

\Prf Write  $A_i = M_i^+ - M_i$ and $B_i = N_i^+ - N_i$, and let
$\theta,\theta_0$ be the functions such that $\fX=\fX(\theta)$ and
$\fX_0=\fX(\theta_0)$. 

If $x \in \H_1$, let $N$ be the smallest subspace of $\N$
containing $x$.
Then $\ol{\fX x} = \theta(N)$.
Thus if $N \ne N_i^+$ for some $i \in \I$, we also have
$\ol{\S x} = \theta(N)$.
When $N = N_i^+$ for some $i \in \I$, 
$\ol{\S x} = \theta(N_i) + [X_i B_ix]$ where $[y]$ denotes the
projection onto $\bC y$.

Suppose that $T \in \B(\H_1,\H_2)$ is given.
Then 
\begin{align*}
  \beta_{\S}(T) &= \sup_{\|x\|=1} \| P_{\S x}^\perp Tx\| \\
  &= \max\big\{ \sup_{N \in \N} \| P_{\theta(N)}^\perp T P_N \|,\ 
                \sup_{i\in \I}\sup_{\substack{x\in N_i^+\\ \|x\|=1}}
                \| (P_{\theta(N_i)}^\perp \!-\! [X_iB_ix] )Tx \|\, \big\} \\
 &= \max\big\{ \beta_{\fX}(T),\
 \sup_{i\in \I} \beta_{\bC X_i}(P_{\theta(N_i)}^\perp T|_{N_i^+} ) \big\} 
\end{align*}

The first important observation is that 
$P_{\theta(N_i)}^\perp \fX P_{N_i^+} = \bC X_i$ is a 1-hyper\-reflexive 
subspace of $\B(N_i^+, \theta(N_i)^\perp)$.
Thus if $T \in \B(\H_1,\H_2)$, there is a multiple $t_i X_i$ so that 
\[ 
 \| P_{\theta(N_i)}^\perp T|_{N_i^+} - t_i X_i \| = 
 \beta_{\bC X_i}(P_{\theta(N_i)}^\perp T|_{N_i^+} ).
\]

Clearly, $\|t_iX_i\|$ is uniformly bounded. Let $S \in \S$ be the
diagonal element given by $SB_i = t_i X_i$.  Replace $T$ by 
$T' = T-S$.  Observe that
\begin{align*}
 \beta_\S(T') &= \beta_\S(T) =
 \max\big\{ \beta_{\fX}(T),\
 \sup_{i\in \I} \beta_{\bC X_i}(P_{\theta(N_i)}^\perp T|_{N_i^+} ) \big\} \\
 &= \max\big\{ \beta_{\fX}(T'),\
 \sup_{i\in \I} \| P_{\theta(N_i)}^\perp T' P_{N_i^+} \| \big\} \\
 &= \beta_{\fX_0}(T') .
\end{align*}

Since $\fX_0$ is 1-hyper\-reflexive, there is an element $X_0 \in
\fX_0$ so that $\| T-(S+X_0) \| = \| T' - X_0 \| = \beta_{\fX_0}(T')$.
This is the desired approximant, showing that $\S$ is
1-hyper\-reflexive.
\qed

The second construction is more elementary.

\begin{lem}\label{diag_const}
Let $\P = \{ P_i : i \in \I \}$ and  $\Q = \{ Q_i : i \in \I \}$ be
partitions of the identity of $\H_1$ and $\H_2$ respectively.
For each $i \in \I$, let $\fX_i$ be a $1$-hyperreflexive subspace of
$Q_i \B(\H) P_i$.
Then the subspace
\[
 \fX := \{ T \in \B(\H) : Q_iTP_i \in \fX_i \qforal i \in \I \} 
\]
is $1$-hyperreflexive.
\end{lem}

\Prf This is straightforward.   Observe that if $x = P_i x$ for some
$i \in \I$ then $\ol{\fX x} = Q_i^\perp \H_2 + \ol{\fX_i x}$;     
while otherwise $\fX x = \H_2$.
Suppose that $T \in \B(\H_1,\H_2)$ and set $T_i = Q_iT|_{P_i\H_1}$.
It is easy to see that 
$\beta_{\fX}(T) = \sup_{i\in \I} \beta_{\fX_i}(T_i)$. 
For each $i \in \I$, there is an $X_i \in \fX_i$ so that
\[ \| T_i - X_i \| = \beta_{\fX_i}(T_i) . \]
Then $X = T + \sum_{i\in \I} Q_i (X_i - T_i) P_i$ lies in $\fX$
and 
\[ \|T-X\| = \| \diag(T_i - X_i) \| = \beta_{\fX}(T) .\]
Thus $\fX$ is 1-hyper\-reflexive.
\qed 

Next, we need a simple way to recognize a nest bimodule.

\begin{prop}\label{P:nest}
A subspace is a nest bimodule if and only if
it is reflexive and the collection of subspaces 
$\{ \ol{\S x} : x \in \H \}$ is totally ordered.
\end{prop}

\Prf If $\S$ is a nest bimodule, then 
$\ol{\S x} = \ol{\T(\M) \S \ol{\T(\N) x}}$.
Now $\ol{\T(\N) x}$ is a subspace $N \in \N$, which is a nested
collection.  Thus its image under $\S$ is also nested.
It follows from the Erdos--Power Theorem \cite{EP} that $\S$ is reflexive.
The map $\theta$ is given by $\theta(N) = \ol{\S N}$.

Conversely, suppose that $\S$ is reflexive and 
$\{ \ol{\S x} : x \in \H \}$ is totally ordered.
Let $\M = \{ \ol{\S x} : x \in \H \} \cup \H$.
For any $S \in \S$ and $T \in \T(\M)$ and $x \in \H$,
$TSx \in T \ol{\S x} \subset \ol{\S x}$.
As $\S$ is reflexive, $\S = \T(\M)\S$.

Now observe that $\S^* = \S^* \T(\M)^* = \S^* \T(\M^\perp)$
is a right nest module.  So the argument of the first paragraph 
shows that the ranges $\ol{\S^*x}$ are totally ordered.
Define $\N^\perp = \{ \ol{\S^* x} : x \in \H \} \cup \H$.
As in the second paragraph, $\S^* = \T(\N^\perp) \S^*$.
Hence $\S = \S \T(\N) =  \T(\M) \S \T(\N)$ is a nest bimodule.
\qed

Here is an easy general condition for $1$-hyperreflexivity.

\begin{prop}\label{P:dist1}
A \wot-closed subspace $\S$ is $1$-hyperreflexive if and only if: 
for every $T\in\B(\H)$ with $\dist(T,\S) = \|T\|$,
there is a sequence of unit vectors $x_n \in \H$ with 
$\dlim_{n\to\infty} \|Tx_n\| = \|T\|$ such that
$\dlim_{n\to\infty} \|P_{\S x_n} Tx_n \| = 0$.
\end{prop}

\Prf By hypothesis, for any $\ep>0$, there is a unit vector $x$
so that $\| P_{\S x}^\perp Tx \| > \|T\| - \ep$.
Therefore $\|Tx\| > \|T\| - \ep$ and 
\begin{align*}
 \| P_{\S x} Tx \|^2 &\le \|Tx\|^2 - \| P_{\S x}^\perp Tx \|^2 \\
 &\le \|T\|^2 - (\|T\| - \ep)^2 < 2 \|T\| \ep .
\end{align*}
The converse is even easier.
\qed

In the case of Magajna's Theorem \cite{Mg}, where $\S = \bC A$, the
condition becomes: $\dlim_{n\to\infty} \|Tx_n\| = \|T\|$ and
$\dlim_{n\to\infty} \ip{Tx_n, Ax_n} = 0$.
He defines a set $W_A(T)$ to be the set of scalars $\lambda \in \bC$
for which there are unit vectors $x_n$ with 
$\dlim_{n\to\infty} \|Tx_n\| = \|T\|$ and
$\dlim_{n\to\infty} \ip{Tx_n, Ax_n} = \lambda$.
The proof proceeds by showing that this set is convex; and if it does
not contain $0$, then a multiple of $A$ may be subtracted from $T$ to
reduce its norm.
There is no obvious way to define such a set for a higher dimensional
algebra that will accomplish the same thing.

\begin{thm}\label{compress}
Let $\S$ be a subspace of $\B(\H_1,\H_2)$.
Suppose that there are projections $P$ and $Q$, at least one of which
is finite rank, so that the compression $Q\S P$ is not $1$-hyperreflexive.
Then neither is $\S$.
\end{thm}

\Prf We may assume that $P$ is finite rank (for if it were $Q$,
we could consider $\S^*$ instead).
Select an element $T\in \B(P\H_1,Q\H_2)$ so that 
\[
  1 = \|T\| = \dist(T, Q\S P) >
 \sup_{\|x\|=1,\,x\in P\H_1} \| P_{Q\S Px}^\perp T x \| =: \beta .
\]
Consider $QTP$ as an element of $\B(\H_1,\H_2)$.  Then 
\[
 \|QTP\| \ge \dist(QTP,\S) \ge \dist(QTP, Q\S P) = \|T\| = \|QTP\|
.\] 
So $\dist(QTP,\S) = 1$.

If $\S$ were 1-hyper\-reflexive, Proposition~\ref{P:dist1}
provides a sequence of unit vectors $x_n \in \H$ so that 
\[
 \dlim_{n\to\infty} \|QTPx_n\| = \|QTP\| \qand
 \dlim_{n\to\infty} \|P_{\S x_n} QTPx_n \| = 0 .
\]
In particular,  
$\dlim_{n\to\infty} \|P^\perp x_n\| = 0$.

Since $\rank(P)$ is finite, there is a subsequence (which we relabel
as $x_n$) so that $x = \dlim_{n\to\infty} x_n$ exists
(and lies in $P\H_1$).
Clearly 
\[
 \| QTP x \| = \lim_{n\to\infty} \|QTP x_n\| = \|QTP\| .
\]
The projections $P_{\S x_n}$ need not converge to $P_{\S x}$;
but there is a lower semicontinuity: if $y \in \ol{\S x}$,
then $\dlim_{n\to\infty} P_{\S x_n} y = y$.
Thus 
\[
 \lim_{n\to\infty} P_{\S x} P_{\S x_n} =
 \lim_{n\to\infty} P_{\S x_n} P_{\S x} = P_{\S x} .
\]
Consequently,
\[
 \| P_{\S x} QTPx\| \le \lim_{n\to\infty} \|P_{\S x_n} QTPx_n \| = 0 .
\]
So 
\[
  \| P_{\S x}^\perp QTP x \| = \| QTPx \| = \|QTP\| .
\]
Therefore $QTPx$ is orthogonal to $\S x$.  As it is obviously 
orthogonal to $Q^\perp \S x$, it is also orthogonal to $Q\S x$.
Hence $\| P_{Q\S x}^\perp Tx\| = \|T\|$ contrary to our hypothesis.
This contradiction establishes the result.
\qed

\begin{rem}
There is no straightforward way to quantify this.
For example, if one takes $\A_n$ to be the algebra of $2 \times 2$ matrices
of the form $\begin{bmatrix}a&n(a\!-\!b)\\0&b\end{bmatrix}$, then it is easy
to check using the matrix $T =
\begin{bmatrix}1&\phantom{-}0\\0&-1\end{bmatrix}$ that 
$\kappa_{\A_n} \ge \dfrac{n^2+1}{2n}$.  
However the infinite inflation $\A_n^{(\infty)} = \A_n \otimes \bC I$
always has distance constant at most $3$ by Bercovici's Theorem \cite{Ber}.
\vspace{.2ex}
Indeed this even holds for the algebra of matrices of the form
$\begin{bmatrix}a&b\\0&a\end{bmatrix}$, which is not even reflexive.
So the compression even to a direct summand of $\A_n^{(\infty)}$ can yield
an arbitrarily large distance constant, or none at all, while the distance
constant for the algebra remains bounded.

Another example of the difficulty in quantifying this can be obtained as
follows. Let $\S_n$ be the subspace of $2\times4$ matrices of the form
$\begin{bmatrix}A & -nA\end{bmatrix}$ for $A \in \A$, a hyperreflexive
subspace. Choose an operator $T \in \fM_2$ so that
$\|T\| = 1 = \dist(T,\A)$ while $\beta_\A(T) = 1/\kappa_\A$.
It would be natural to try $T' = \begin{bmatrix}T&0\end{bmatrix}$ as
a test case for the distance constant for $\S_n$.
Clearly $\|T'\| = 1 = \dist(T',\S_n)$.
However pick unit vectors $x$ and $y$ so that $Tx=y$.
Observe that $x' = \begin{bmatrix}nx/\sqrt{n^2+1}\\
x/\sqrt{n^2+1}\end{bmatrix}$ belongs to $\ker \S_n$.  
Thus $\beta_{\S_n}(T') \ge y^* T' x' = \dfrac n {\sqrt{n^2+1}}$.
So the \textit{proof} of Theorem~\ref{compress} does not reveal much about
the distance constant of $\S_n$. 
Nevertheless, in this example, one can show that  
$\kappa_{\S_n} \ge \kappa_\A$.
\end{rem}

The following result shows that, for bimodules over masas,
the 1-hyper\-re\-flex\-ive ones are obtained using Lemma~\ref{diag_const}
where the diagonal entries are nest bimodules.
(The only 1-dimensional $\fC$--$\fD$ masa bimodule has the form
$Q \B(\H)P$ where $Q$ and $P$ are one dimensional projections.  
This is a nest bimodule. 
So type S2 reduces to type S1 in this case.)

\begin{thm}\label{masa modules} 
Suppose that $\fX$ is a $1$-hyperreflexive subspace such that 
the families of projections $\{ P_{\fX x} : x \in \H_1 \}$
and $\{ P_{\fX^* y} : y \in \H_2 \}$ are both commutative.
Then there are abelian von Neumann algebras $\fC$ and $\fD$ in $\B(\H_2)$
and $\B(\H_1)$ respectively so that $\fX$ is a $\fC'$--$\fD'$ bimodule. 
Moreover, there are two collections of pairwise orthogonal projections 
\[
 \C = \{ C_j : j \in \J \} \subset \fC \qand 
 \D = \{ D_j : j \in \J \} \subset \fD 
\]
and nest bimodules $\fX_j \subset C_j \B(\H_1,\H_2) D_j$ so that 
\[
 \fX = \{ T \in \B(\H) : C_jTD_j \in \fX_j \qforal j \in \J \} .
\]
\end{thm}

\Prf Since the projections $P_{\fX x}$ commute, there is an abelian von
Neumann algebra $\fC$ in $\B(\H_2)$ containing all of them.  
Observe that if $X \in \fX$, $C \in \fC'$ and $x \in \H_1$, 
then 
\[
 CXx \in C P_{\fX x}\H_2 = P_{\fX x} C \H_2 \subset \ol{\fX x} .
\]
As $\fX$ is reflexive, $\fC' \fX = \fX$.
Similarly there is an abelian von Neumann algebra $\fD$ in $\B(\H_1)$
containing  $\{ P_{\fX^* y} : y \in \H_2 \}$, and $\fX = \fX \fD'$.

For each non-zero vector $x \in \H_1$ such that $\ol{\fX x} \ne \H_2$,
let $C_x$ be the smallest projection in $\fC$ such that 
\begin{enumerate}
\item $C_x (\fX x)^\perp = (\fX x)^\perp$, and 
\item for all $y \in \H_1$, either $C_x (\fX y)^\perp = 0$ or 
 $C_x (\fX y)^\perp = (\fX y)^\perp$.
\end{enumerate}
There is such a smallest projection because the product of any two
projections with this property also has the property; and so does the
(decreasing) limit of any sequence of such projections.

Let $D_x$ be the projection onto 
\[ \spn \{y \in \H_1 : C_x (\fX y)^\perp \ne 0 \} . \]
Then $D_x \in \fD$.
Indeed, if $C_x (\fX y)^\perp \ne 0$ and $z = Dy$ for $D \in \fD'$,
then $\fX z \subset \fX y$ and hence 
$(\fX z)^\perp \supset (\fX y)^\perp \ne 0$.
Hence $C_x (\fX z)^\perp \ne 0$ and $z = D_x z$.
Thus the range of $D_x$ is $\fD'$--invariant; 
so $D_x$ lies in $\fD'' = \fD$.

Observe that for vectors $x$ and $y$, either $C_x=C_y$
or $C_xC_y=0$.
Indeed, this follows from the minimality of $C_x$ and $C_y$.
For if $C_xC_y\ne 0$, then either $C_y(\fX x)^\perp = 0$ and
thus $C_x C_y^\perp$ will be a smaller projection satisfying the two
conditions, a contradiction; or $C_y(\fX x)^\perp = (\fX x)^\perp$
so that $C_xC_y$ is such a projection.  
This latter condition is not contradictory only if $C_x \le C_y$.
But by symmetry, we also obtain $C_y \le C_x$, whence equality.

Since $D_x$ is a function only of $C_x$, not $x$ itself, we obtain that
$D_x=D_y$ when $C_x=C_y$.  
If $C_x C_y=0$, then $D_xD_y=0$ also.
To see this, it suffices to show that if 
\[
 C_x (\fX u)^\perp = (\fX u)^\perp \qand
 C_y (\fX v)^\perp = (\fX v)^\perp ,
\]
then $\ip{u,v}=0$. But if this inner product is non-zero,
the two $\fD'$-modules $\ol{\fD' u}$ and $\ol{\fD' v}$ have non-trivial
intersection, say containing a non-zero vector $w$.
The set of vectors $z$ with $C_x (\fX z)^\perp = (\fX z)^\perp$
was shown to be invariant under multiplication by $\fD'$ and is clearly
norm closed.  Thus $w$ is in this set, so that
$C_x (\fX w)^\perp = (\fX w)^\perp$.
Similarly $C_y (\fX w)^\perp = (\fX w)^\perp$; and thus $C_xC_y \ne 0$.

Therefore the collection of all projections $D_x$ may be
enumerated as a family $\{D_j : j \in \J\}$ of pairwise orthogonal 
projections with corresponding projections $\{C_j : j \in \J \}$.
If $\ol{\fX x} \ne \H_2$, then there is an $j \in \J$ so that $D_j =
D_x$ and $C_j = C_x$.  
Hence $\ol{\fX x}$ contains $C_j^\perp \H_2$.
Define spaces $\fX_j = C_j \fX D_j$ for $j \in \J$.
For every vector $0 \ne x \in D_j\H_1$, $\ol{\fX x}$ contains $C_j^\perp
\H_2$. Since $\fX$ is reflexive and is a $\fC'$--$\fD'$ module, 
it must contain $C_j^\perp \B(\H_1,\H_2) D_j$ and $\fX_j$.
It follows that
\[ \fX = \{ X \in \B(\H) : C_j X D_j \in \fX_j,\ j \in \J \} .\]

The last step is to show that each $\fX_j$ is a nest bimodule.
Fix $j \in \J$ and work in $\B(D_j\H_1, C_j\H_2)$.
Let $\fC_j$ and $\fD_j$ be the restrictions of the abelian von Neumann
algebras. By Proposition~\ref{P:nest}, it suffices to show that the 
projections $Q_x = P_{\ol{\fX_j x}}$ in $\fC_j$ are nested for 
$x \in D_j\H_1$.

Suppose to the contrary that there are vectors $x$ and $y$
for which this fails.  
So $Q_x Q_y^\perp \ne 0$ and $Q_x^\perp Q_y \ne 0$.
Now $Q_x \H_2 = \ol{ \fX_j \fD_j x}$ depends only on the projection 
$P_x$ onto $\ol{\fD'_j x}$ in $\fD_j$.
We first show that we may suppose $P_x P_y = 0$.
Indeed, if this is not the case, we may choose vectors $x_1, y_1,z$ so
that $x_1$ has support $P_xP_y^\perp$, $y_1$ has support $P_y
P_x^\perp$ and $z$ has support $P_x P_y$.
Observe that 
\[
  Q_x = Q_{x_1} \vee Q_z \qand Q_y = Q_{y_1} \vee Q_z ,
\]
whence 
\setcounter{equation}{2}
\begin{align}\label{E3}
  Q_x^\perp = Q_{x_1}^\perp Q_z^\perp \qand
  Q_y^\perp = Q_{y_1}^\perp Q_z^\perp .
\end{align}
So 
\[
 0 \ne Q_x Q_y^\perp = (Q_{x_1} \vee Q_z) Q_{y_1}^\perp Q_z^\perp
  = Q_{x_1} Q_{y_1}^\perp Q_z^\perp .
\]
So $Q_{x_1} Q_{y_1}^\perp \ne 0$ and likewise 
$Q_{x_1}^\perp Q_{y_1} \ne 0$.
This reduces us to vectors with disjoint support projections.

Consider the restriction of $\fX_j$ to 
$(P_{x_1} \oplus P_{y_1}) \H_1$ and compress the range to 
\[
 Q_{x_1}^\perp \vee Q_{y_1}^\perp = 
 Q_{x_1} Q_{y_1}^\perp + Q_{x_1}^\perp Q_{y_1}
 + Q_{x_1}^\perp Q_{y_1}^\perp .
\]
Then $(Q_{x_1}^\perp \vee Q_{y_1}^\perp)\fX_i (P_{x_1} + P_{y_1})$ 
has the form
\[
 \begin{bmatrix} * & 0\\0 & *\\0 & 0\end{bmatrix}
\]
where the $*$ entries are non-zero and independent.
If $Q_{x_1}^\perp Q_{y_1}^\perp \ne 0$, we may choose vectors in each
subspace and compress.
This will be a $3\times2$ masa bimodule of the same form, which has
distance constant $\sqrt{9/8}$ (see~\cite{DO}).
So by Theorem~\ref{compress}, $\fX_i$ would not be 1-hyper\-reflexive.
Thus $Q_{x_1}^\perp Q_{y_1}^\perp = 0$.

To recap, this shows that if ${x_1}$ and ${y_1}$ have disjoint supports,
then either $Q_{x_1}^\perp$ and $Q_{y_1}^\perp$ are comparable or they 
are orthogonal.
Equation~(\ref{E3}) shows that this remains true if we drop the
condition on disjoint supports.

Fix ${x_1}$ as above with $Q_{x_1}^\perp \ne 0$, and consider the span
$Q$ of all projections $Q_z^\perp$ which are comparable to $Q_{x_1}^\perp$.
For each vector $y$ such that $Q_y^\perp Q \ne 0$, one has 
$Q_y^\perp Q_z^\perp \ne 0$ for some $Q_z^\perp$ which is comparable to
$Q_{x_1}^\perp$. So either $Q_y^\perp \le Q_z^\perp \le Q$; or
$Q_y^\perp \ge Q_z^\perp$ and thus $Q_y^\perp Q_{x_1}^\perp \ne 0$.
Hence $Q_y^\perp$ is comparable to $Q_{x_1}^\perp$ and thus is less than
$Q$. So $Q$ is a projection satisfying conditions (1) and (2), for which
$C_j$ is the minimal choice.  
But we are working in $\B(C_j\H_1, D_j\H_2)$ so $Q \le C_j$ is automatic. 
Thus $Q=C_j$.

Return to $Q_{y_1}^\perp$, which is orthogonal to $Q_{x_1}^\perp$.
By the previous paragraph, there is a vector $z_1$ so that
$Q_{z_1}^\perp$ is comparable to $Q_{x_1}^\perp$ and
$Q_{y_1}^\perp Q_{z_1}^\perp \ne 0$.  
We deduce that $Q_{z_1}^\perp > Q_{y_1}^\perp$.
Observe that, since $Q_{x_1}Q_{y_1}^\perp \ne 0$ and 
$Q_{z_1}Q_{y_1}^\perp = 0$, if we replace $x_1$ by 
$x_2 = P_{z_1}^\perp x_1$ then $Q_{x_2}Q_{y_1}^\perp \ne 0$.
Since $Q_{x_2} \le Q_{x_1}$, $Q_{x_2}^\perp Q_{y_1} \ne 0$ too.
Similarly we can replace $y_1$ by $y_2 = P_{z_1}^\perp y_1$
and still maintain the fact that 
\[ Q_{x_2}Q_{y_2}^\perp \ne 0 \qand Q_{x_2}^\perp Q_{y_2} \ne 0 .\]
Moreover since 
$
Q_{z_1}^\perp Q_{x_2}^\perp \ne 0$ and $Q_{z_1}^\perp Q_{y_2}^\perp \ne 0,
$
it follows again that $Q_{z_1}^\perp \ge Q_{x_2}^\perp + Q_{y_2}^\perp$.
But now $x_2$, $y_2$ and $z_1$ have disjoint supports.

Pick vectors $u \in Q_{x_2}^\perp Q_{y_2}\H_2$ and 
$v \in Q_{x_2} Q_{y_2}^\perp \H_2$; and consider the compression 
$\S = P_{\spn\{u,v\}} \fX_j P_{\spn\{ x_2, y_2, z_1 \}}$.
This has the form $\begin{bmatrix}0 & * & 0\\ * & 0 & 0\end{bmatrix}$
with respect to this decomposition and the two non-zero entries are
not dependent.  
Thus this has distance constant $\sqrt{9/8}$.
By Theorem~\ref{compress}, $\fX_j$ is not 1-hyper\-reflexive.
This contradiction establishes the fact that $\fX_j$ is a nest bimodule.
\qed

For future use, we record one fact that is a consequence of the proof.
The $\fC_i$--$\fD_i$ module $\fX_i$ not only has the property that
the range projections $\{ P_{\fX_i x} \}$ and $\{ P_{\fX^* u} \}$
are commutative.  It also has the minimality hypothesis that the
projections
$C_x$ satisfying (1) and (2) above are $I$ or $0$.  
When $\fX_i$ was not a nest bimodule, this allowed us to find
compressions which were $2\times3$ or $3\times2$ submodules which are
evidently not 1-hyperreflexive.

\begin{cor}\label{C:masa}
Assume that $\fX$ is a $\fC'$--$\fD'$ bimodule for two abelian von 
Neumann algebras $\fC$ and $\fD$, 
Assume also that no proper projection $C \in \fC$ has the  property that
$C (\fX x)^\perp$ is either $0$ or $(\fX x)^\perp$
for each $x \in \H_2$. 
If $\fX$ is not a nest bimodule, then either
\begin{enumerate}
\item[(a)] there are orthogonal projections $D_1, D_2 \in \fD$ and $C_1,
C_2, C_3 \in \fC$ so that $(C_1+C_2+C_3) \fX (D_1+D_2)$ has the form
$
 \begin{bmatrix} * & 0\\0 & *\\0 & 0\end{bmatrix}
$
where the $*$ entries are non-zero and independent, or

\item[(b)] there are orthogonal projections $D_1, D_2, D_3 \in \fD$ and
$C_1, C_2 \in \fC$ so that $(C_1+C_2) \fX (D_1+D_2+D_3)$ has the form
$
 \begin{bmatrix} * & 0 & 0\\0 & * & 0\end{bmatrix}
$
where the $*$ entries are non-zero and independent.
\end{enumerate}
\end{cor}

\section{1-Hyperreflexive subspaces in low dimensions}

In order to conveniently eliminate subspaces as failing to be
1-hyper\-re\-flex\-ive, we need some low dimensional examples.
In this section, we characterize the subspaces of $2\times2$ and $2\times3$
matrices which are 1-hyper\-reflexive.

\begin{thm}\label{T:22case}
Let $\S$ be a $1$-hyper\-re\-flex\-ive subspace of $\fM_2$.
Then it is one of the following:
\begin{enumerate}
\item $\Dim \S$ is $0$,$1$ or $4$.
\item $\Dim \S = 3$ and there are unit vectors $x$ and $y$
so that\\ $\S = \{ T \in \fM_2 : \ip{Tx,y}=0 \}$.
\item $\Dim \S = 2$ and
\begin{enumerate}
\item $\Dim \ran \S = 1$, so that $\S = Q \fM_2$ for some projection $Q$.
\item $\Dim \ker \S = 1$, so that $\S = \fM_2 P$ for some projection $P$.
\item there are unitaries $U$ and $V$ so that $\S = U \fD_2 V$,\\
where $\fD_2$ is a masa.
\end{enumerate}
\end{enumerate}
Cases $(1)_0$, $(1)_4$, $(2)$, $(3\mathrm{a})$ and $(3\mathrm{b})$
are nest bimodules; and case $(1)_1$ is $1$-dimensional.
Case $(3\mathrm{c})$ is type $\mathrm{S3}$.
\end{thm}

\pagebreak[3]
\Prf The cases of $\Dim\S = 0$ or $4$ are trivial, and $\Dim\S = 1$ 
is Magajna's Theorem.

If $\Dim\S=3$ and there is a vector $x$ so that $\Dim\S x = 1$,
then by a dimension count, one concludes that 
$\S = \{ T \in \fM_2 : \ip{Tx,y}=0 \}$ where $y$ is chosen orthogonal
to $\S x$.  Evidently $\S u$ belongs to $\{ 0, \S x, \H \}$ depending
on whether $u=0$, $u \in \bC^* x$, or not, respectively.
By Proposition~\ref{P:nest}, $\S$ is 1-hyper\-re\-flex\-ive.

On the other hand, if $\Dim\S x = 2$ for all $x \ne 0$, then
$\Ref(\S) = \fM_2$.  So $\S$ is not even reflexive.

Now consider $\Dim \S = 2$.
Cases (3a) and (3b) are evidently 1-hyper\-re\-flex\-ive.
They are both nest bimodules.
So we assume that $\S$ has no proper kernel or cokernel.

As in the 3-dimensional case, there must be a vector $x_1$
so that $\S x_1 = \bC y_1$ is 1-dimensional.
If all vectors $x \not\in \bC x_1$ had 2-dimensional range under $\S$,
the 3-dimensional case again shows that $\Ref(\S)$ would be
a 3-dimensional nest bimodule.  So there is a second vector $x_2$
independent of $x_1$ so that $\S x_2 = \bC y_2$.
If $y_2$ is a multiple of $y_1$, it follows that $\S \H = \bC y_1$,
which is case (a).  Thus we have $y_1$ and $y_2$ independent.

Next observe that the functionals $\phi_i(S) = \ip{Sx_i,y_i}$
must be independent.  For otherwise, $\S$ would be 1-dimensional.
Consider $\S$ with respect to an orthogonal basis $x_1, x_1'$
for the domain and $y_1, y_1'$ for the range.  
Then $\phi_1'(S) = \ip{Sx_1',y_1'}$ is easily seen to be 
independent of $\phi_1$.
In this basis, we now have independent entries in the $1,1$
and $2,2$ positions and $0$ in the $2,1$ position.
The $1,2$ entry must be a linear combination of the other two.
Hence there are scalars $r$ and $s$ so that
\[ \S = \begin{bmatrix}a&ar+bs\\0&b\end{bmatrix}.\]

We may further simplify this to the case of $r\ge0$ and $s\ge0$
as follows.  
Write $r=|r|e^{i\rho}$ and $s=|s|e^{i\sigma}$.
Then
\[
 \begin{bmatrix}1&0\\0&e^{i\sigma}\end{bmatrix}
 \begin{bmatrix}a&ar+bs\\0&b\end{bmatrix}
 \begin{bmatrix}1&0\\0&e^{-i\rho}\end{bmatrix}
= \begin{bmatrix}a&a|r|+be^{i(\sigma-\rho)}|s|\\
  0&be^{i(\sigma-\rho)}\end{bmatrix}.
\]
In this basis, we see that $x_2 =
\begin{spmatrix}-r\\ \phantom{-}1\end{spmatrix}$ and
$y_2 = \begin{spmatrix}s\\1\end{spmatrix}$.

Observe that $\S x = \H$ except when $x$ is
a multiple of either $x_1$ or $x_2$.
Hence if we select a unit vector $y_2'$ orthogonal to $y_2$,
\begin{align*}
 \beta_\S(T) &= \max_{i=1,2}  \| P_{\S x_i}^\perp Tx_i\| 
 = \max_{i=1,2}  |\ip{Tx_i, y_i'}| 
\end{align*}

The proof will be complete once we show that
1-hyper\-re\-flex\-ivity implies that $r=s=0$.
Suppose that we have a unitary operator $U$
such that $\dist(U,\S)=1$.  
We will have $\beta_\S(U)<1$ if and only if
$\ip{Ux_i,y_i} \ne 0$ for $i=1,2$.

Define a unitary $U = \begin{bmatrix}\alpha & -\beta\\ \beta&\alpha\end{bmatrix}$
where $\alpha = \sin \theta$ and $\beta = \cos\theta$ satisfy
\[
 0 < \alpha < (r+s)\beta \qand \ip{Ux_2,y_2} \ne 0 .
\]
Then $\beta_\S(U)<1$.  We claim that $\dist(U,\S)=1$. 

Suppose to the contrary that there are scalars $a$ and $b$ so that
\[
 \left\| \begin{bmatrix}\alpha-a & -\beta-ra-sb\\ \beta&\alpha-b\end{bmatrix} \right\| < 1 .
\]
We may suppose that $a$ and $b$ are real since the complex
conjugate will have the same norm, and one can average to
replace $a$ and $b$ by their real parts while decreasing the norm.
Clearly $a$ and $b$ are strictly positive, for otherwise either the
first column or second row will have norm at least one.
The first row will have norm less than 1, and so
\[
 (\alpha-a)^2 + (-\beta-ra-sb)^2 < 1
\]
whence
\[
 2a\alpha > 2\beta(ra+sb) + (ra+sb)^2 + a^2 .
\]
Similarly the second column leads to the inequality
\[
 2b\alpha > 2\beta(ra+sb) + (ra+sb)^2 + b^2 .
\]
Multiply the first by $r$ and the second by $s$, 
add them and divide by $2(ra+sb)$ to obtain
\[
 \alpha
  > (r+s)\beta + \frac12 (ra+sb)(r+s) + \frac{ra^2+sb^2}{2(ra+sb)}
 > (r+s)\beta
\]
This is a contradiction, which establishes our claim.
Thus 1-hyper\-re\-flex\-ivity shows that $r=s=0$, which is case (c).
\qed

\begin{thm}\label{T:23case}
Let $\S$ be a $1$-hyper\-re\-flex\-ive subspace of $\fM_{2,3}$.
Then it is one of the following:
\begin{enumerate}
\item $\S$ is a nest bimodule.

\item $\Dim \S = 1$.

\item $\Dim\S=4$ and there are orthonormal bases
so that it consists of all matrices of the form
$\begin{bmatrix}a&0&c\\0&b&d\end{bmatrix}$.
This is type $\mathrm{S3}$.

\item $\Dim\S = 3$ and there are orthonormal bases
so that it consists of all matrices of the form
$\begin{bmatrix}a&0&0\\0&b&c\end{bmatrix}$.
This is type $\mathrm{S3}$.

\item $\Dim\S = 3$ and there are orthonormal bases
so that it consists of all matrices of the form
$\begin{bmatrix}aT & \begin{matrix}b\\c\end{matrix}\end{bmatrix}$
for some $T \in \fM_2$ of rank $2$.
This is an instance of the construction of Lemma~$\ref{diag_const}$. 

\end{enumerate}
\end{thm}

\Prf Suppose that $\S$ is reflexive and $\Dim\S \ge 2$ but $\S$ is not
a nest bimodule.
By Proposition~\ref{P:nest}, 
the subspaces  $\{\S x : x \in \H\}$ are not totally ordered.  
Then there are two vectors $x_1$ and $x_2$ so that
$\S x_1 = \bC y_1$ and $\S x_2 = \bC y_2$, where
$y_1$ and $y_2$ are independent.
So $\Dim\S \le 4$.

First assume that the two functionals $\phi_i(S) = \ip{Sx_i,y_i}$
are independent on $\S$.
By Theorem~\ref{compress}, the compression of the domain
to $\spn\{x_1,x_2\}$ must be 1-hyper\-re\-flex\-ive.
Thus by Theorem~\ref{T:22case}, this forces $x_1$ and $x_2$
to be orthogonal and likewise $y_1$ and $y_2$ are orthogonal.
Thus $\S$ has the form 
$\begin{bmatrix}a&0&?\\0&b&?\end{bmatrix}$.

When $\phi_1$ and $\phi_2$ are dependent, let $P$ be the
projection onto the $\spn\{x_1,x_2\}$.
Then $\S P$ is one dimensional, so equals $\bC T$ for
some $2\times 2$ matrix $T$ of rank 2.
In this case, $\Dim\S \le 3$ and has the form
$\begin{bmatrix}aT & \begin{matrix}?\\?\end{matrix}\end{bmatrix}$.

When $\Dim\S = 4$, the functionals are indeed independent.
Therefore this puts $\S$ into the predicted form
$\begin{bmatrix}a&0&c\\0&b&d\end{bmatrix}$.

When $\Dim\S = 3$ and $\phi_1$ and $\phi_2$ are dependent,
the unknowns are independent variables and we have the form
$\begin{bmatrix}aT & \begin{matrix}b\\c\end{matrix}\end{bmatrix}$.

When $\Dim\S = 3$ and $\phi_1$ and $\phi_2$ are independent, 
at least one of the coefficients marked ?\
will be independent of $a$ and $b$.  By symmetry, we may suppose
the form $\begin{bmatrix}a&0&L(a,b,c)\\0&b&c\end{bmatrix}$
\vspace{.2ex}
where $c$ is independent of $a$ and $b$ and $L(a,b,c) = ra+sb+tc$ is linear.
Compress to $\S P_1$ where $P_1$ is the projection onto the
subspace $\spn\{x_1, sx_2+tx_3\}$. 
This yields
$\S P_1 \simeq \begin{bmatrix}a&t(ra + d) \\0&d\end{bmatrix}$
where $d=sb+tc$.
By Theorem~\ref{T:22case}, we deduce that $t=0$.

Similarly compress by the projection $P_2$ onto
$\spn \big\{ (x_1+x_2)/\sqrt2, x_3 \big\}$ to obtain
$
 \S P_2 \simeq
 \begin{bmatrix} a/\sqrt2 & ra+sb\\
 b/\sqrt2 & c \end{bmatrix}
$.
Again by Theorem~\ref{compress} and Theorem~\ref{T:22case},
this three dimensional space must have a vector 
\[ v = \alpha \frac{x_1+x_2}{\sqrt2} + \beta x_3 \]
with one-dimensional range.
Clearly $\beta \ne 0$.
Thus the second coordinate of $SP_2 v$ is arbitary independent
of the first coordinate because $c$ is arbitrary.
So the first coordinate, namely 
$(\alpha/\sqrt2 +r\beta)a + (s\beta)b$,
needs to be zero for all $a$ and $b$.
This forces $s=0$.

So $\S = \begin{bmatrix}a&0&ra\\0&b&c\end{bmatrix}$.
Now a third application of Theorem~\ref{compress} and
Theorem~\ref{T:22case}, compressing to the subspace
$\spn\{x_1,x_3\}$, shows that $r=0$. 
This puts $\S$ in the desired form.
\smallbreak

Now consider the case of $\Dim\S=2$ with $\phi_1$ and $\phi_2$
dependent. 
There is a norm one element $S \in \S$ such that $SP=0$.
Thus we can choose an orthonormal basis $x_1,x_2,x_3$ for the 
domain such that $\spn\{x_1,x_2\} = P\H_1$ and an orthonormal basis
for the range, $y_1,y_2$ so that $S x_3 = y_1$.
Then $\S$ has the form $\begin{bmatrix} aT 
&\begin{matrix} b\\ra \end{matrix}\end{bmatrix}$.
Since $T$ has rank~2, choose a unit vector $x = Px$ so that
$z=Tx = \begin{spmatrix}z_1\\z_2\end{spmatrix}$ is neither collinear nor
orthogonal to $y_1$; i.e.\ $z_1 z_2 \ne 0$. 
\vspace{.4ex}
Compress the domain to $\spn\{x,x_3\}$ via the projection $P_1$. 
Then $\S P_1 = \begin{bmatrix} az_1 & b\\ az_2 & ra \end{bmatrix}$. 
\vspace{.2ex}
If $r\ne0$, this is not reflexive as the subspace ranges:
$0$, $\bC z$ and $\H_2$ are nested.  While if $r=0$, then
$\S P_1$ is diagonalizable but not with orthonormal bases.
Thus by Theorem~\ref{T:22case}, it is not 1-hyper\-re\-flex\-ive.
Hence by Theorem~\ref{compress}, neither is $\S$.

Thus we may suppose that $\phi_1$ and $\phi_2$ are independent.
So by the earlier analysis, $\S$ has the form 
$\begin{bmatrix}a&0&ra+sb\\0&b&ta+ub\end{bmatrix}$.
Compressing the domain to 
\[
 P \H =
 \spn \Big\{ \frac{x_1+x_2}{\sqrt2},
 \alpha x_1 - \alpha x_2 + \beta x_3 \Big\}
\]
where $2 \alpha^2 + \beta^2=1$, will yield
$\begin{bmatrix}a/\sqrt2&(\alpha + r\beta)a + s\beta b\\
b/\sqrt2& t\beta a + (u\beta - \alpha)b\end{bmatrix}$.
This is never 1-hyper\-re\-flex\-ive for all choices of parameters
$\alpha$ and $\beta$.
Indeed, the subspace will have kernel only if $s=t=0$ and   
$\alpha +r\beta = -\alpha + u\beta$.  
In the other cases, one looks for vectors with one dimensional
range.  Generically there are only two such vectors but they are
usually neither parallel nor orthogonal.  
\qed

\section{The Noncommuting case}\label{S:noncomm}

We now have the tools we need to consider 1-hyper\-reflexive
subspaces for which the projections $P_{\S x}$ do not commute.
Once we understand exactly how this can occur, we will be able to complete
the proof of Theorem~\ref{T:1hyp}.

\begin{lem}\label{noncomm1}
Let $\S \subset \B(\H_1,\H_2)$ be a $1$-hyperreflexive subspace.
Suppose that the orthogonal projections 
$Q_x = P_{\S x}$ and $Q_y = P_{\S y}$ do not commute.
Then $Q = (Q_x \vee Q_y) - (Q_x \wedge Q_y)$ has rank two,
and $Q \S |_{\spn\{x,y\}}$ is one dimensional.
\end{lem}

\Prf By Theorem~\ref{compress}, the subspace 
$Q \S |_{\spn\{x,y\}}$ is 1-hyper\-reflexive.
So we may work with this space, so that
$\H_1 = \spn\{x,y\}$, $Q_x \wedge Q_y = 0$
and $Q=Q_x \vee Q_y = I_{\H_2}$.

As $Q_x$ and $Q_y$ do not commute, there is a unit vector
$u = Q_x u$ so that $u \ne v = Q_y u \ne 0$. 
Let $Q_0$ be the projection onto $\spn\{u,v\}$.
Choose a unit vector $u'$ in $\spn\{u,v\}$ orthogonal to $u$.
Also choose an orthonormal basis $\{x,x'\}$ for $\spn\{x,y\}$.
Then there are constants $\gamma$ and $\delta$ so that
$y$ is a non-zero multiple of $\gamma x + x'$
and $v$ is a non-zero multiple of $\delta u + u'$.

Consider the compression $Q_0 \S$.
Since $Q_0 \S x = \bC u$ and $Q_0 \S y = \bC v$, we obtain
that $Q_0 \S \subset \spn\{ ux^*, vy^* \}$.
With respect to the  orthonormal bases $\{x,x'\}$ and $\{u,u'\}$,
$\spn\{ ux^*, vy^* \}$ has the form
\[
 a\begin{bmatrix} 1& 0\\ 0& 0\end{bmatrix} +
 b\begin{bmatrix} \ol{\gamma}\delta& \ol{\gamma}\\ \delta& 1\end{bmatrix} .
\]
This space is 1-hyper\-reflexive.  By Theorem~\ref{T:22case},
it is either one-dim\-en\-sional or the two bases $\{x,y\}$
and $\{u,v\}$ are orthogonal.
As the latter is not the case, this compression is one-dimensional,
say multiples of an operator $T$.

Now suppose that at least one of $Q_x$ or $Q_y$ has rank at least 2,
say $Q_x$.
Then we may select another unit vector $u''$ orthogonal
to $u$, $u'$ and the range of $Q_y$.
Thus $u$ and $Q_x u''$ will be independent, and thus
the compression of $\S x$ to $\spn\{u,u''\}$ is two-dimensional.
Consider the projection $Q_1$ onto $\spn\{u,u',u''\}$.
Then $Q_1 \S$ has the form 
$\begin{bmatrix}aT\\ \begin{matrix}b&0\end{matrix}\end{bmatrix}$.
By Theorem~\ref{T:23case}, this is not 1-hyper\-reflexive.
This contradiction shows that $Q_x$ and $Q_y$ both have rank one.

Therefore $Q_x \vee Q_y$ has rank two, and we have already shown that
$Q \S |_{\spn\{x,y\}}$ is one dimensional.
\qed

\begin{lem}\label{noncomm2}
Let $\S \subset \B(\H_1,\H_2)$ be a $1$-hyperreflexive subspace.
Suppose that the projection $Q_x$ does not commute
with $Q_y$ nor $Q_z$, where $\{x,y,z\}$ are independent.
Then $Q_x \wedge Q_y = Q_x \wedge Q_z.$
Let 
\[ Q = (Q_x \vee Q_y \vee Q_z) - (Q_x \wedge Q_y  \wedge Q_z) .\]
Then the compression $Q\S|_{\spn\{x,y,z\}}$ is one-dimensional.
\end{lem}

\Prf We will work with the compression of the domain of $\S$ to
$\spn\{x,y,z\}$ and the range to $Q\H_2$.
Let $x,x'$ and $u,u'$ be orthonormal bases for $\spn\{x,y\}$
and the range of $R_1 = (Q_x \vee Q_y) - (Q_x \wedge Q_y)$, respectively.
By the previous lemma, $R_1 \S|_{\spn\{x,y\}} = \bC T$
where $T$ is a $2\times2$ matrix of rank two. 
Let $R_0 = Q_x \wedge Q_y$.
For each non-zero vector $\alpha x + \beta x'$, we have
$Q_{\alpha x + \beta x'} = R_0 + [T(\alpha x + \beta x')]$.

Suppose that $R_0 \ne 0$.
Then $0 = R_0 \wedge Q_z$.
However we know from the previous lemma that $Q_x \wedge Q_z$
is codimension one in $Q_x$.
Hence we deduce that $\rank R_0=1$, say spanned by a unit vector $v$
(which is orthogonal to $u,u'$).
Moreover $Q_z$ is a projection of rank 2 onto a subspace that
does not contain $v$.  
But $Q_x \wedge Q_z$ is rank one, so maps onto the span of some
vector $w = \gamma v + \delta Tx$, with $\delta \ne 0$.
If $Q_z$ commutes with~$R_0$, then $R_0Q_z=0$. This would force
$\gamma=0$ and so $Q_z\geq [Tx]$. But then $Q_x=R_0+[Tx]$ would
commute with~$Q_z$. So $Q_z$ does not commute with $R_0$.

Suppose that there were two independent vectors $\alpha_1 x + \beta_1 x'$
and $\alpha_2 x + \beta_2 x'$ so that both $Q_{\alpha_1 x + \beta_1 x'}$
and $Q_{\alpha_2 x + \beta_2 x'}$ commute with $Q_z$.
Then $Q_z$ would also commute with
\[ Q_{\alpha_1 x + \beta_1 x'} \wedge Q_{\alpha_2 x + \beta_2 x'} = R_0 ,\]
contrary to fact.
It follows that that $Q_z$ does not commute with all but at most one of the
projections $Q_{\alpha x + \beta x'}$.
Therefore the previous lemma shows that 
$Q_z \wedge Q_{\alpha x + \beta x'}$ is rank one.
For $\beta \ne 0$, this will not be the vector $w$.  
So the range of $Q_z$ contains a second independent vector in the range of
$R_0+R_1$.  Hence $Q_z \le R_0+R_1$.
As this range does not contain $v$, $Q_1 Q_z$ maps onto the range of $Q_1$.

With respect to the bases $x,x',z$ (which is not orthogonal) and
$u,u',v$ (which is orthonormal), $\S$ has the form
\[
 \begin{bmatrix}aT& \begin{matrix} b\\ c\end{matrix}\\
 \begin{matrix} d & e \end{matrix}& L(b,c)\end{bmatrix} 
\]
where $L(b,c)$ is a linear function.

The discussion above shows that $b$ and $c$ are independent.
It is also the case that they are independent of $a$.
Indeed, the first two rows are 1-hyper\-reflexive.
So it follows from Theorem~\ref{T:23case} as this space must
be three dimensional.
Likewise the $d$ and $e$ are not dependent on $a$ because 
$Q_x$ and $Q_{x'}$ are rank 2.  So by the same reasoning, they
are also independent of each other.

Restrict the domain of $\S$ to the subspace
$\spn\{ \alpha x + \beta x', z \}$
and write $u'' = T(\alpha x + \beta x')$ and $d' = \alpha d + \beta e$. 
Then we obtain
\[
 \begin{bmatrix}a u'' & \begin{matrix} b\\ c\end{matrix}\\
 d' & L(b,c) \end{bmatrix} .
\]
By Theorem~\ref{T:23case}, since this is a 1-hyper\-reflexive space
and is 4 dimensional,
it must be the case that $d'$ is independent of $a,b,c$, and
the functional $L = 0$.  This means that $Q_z = R_1$, a contradiction.

All of this analysis leads to the conclusion that in fact $R_0=0$,
which is to say that $Q_x \wedge Q_y = Q_x \wedge Q_z$.

Therefore each projection $Q_x$, $Q_y$ and $Q_z$ is one dimensional.
The restriction of the domain of $\S$ to $\spn\{ x, z \}$ is one
dimensional. 
Thus for $S \in \S$, $Sz$ is a linear function of $Sx$, as is $Sx'$.
Selecting $S_0$ with $S_0 x \ne 0$, define an
operator $T'=\begin{bmatrix} S_0 x & S_0 x' & S_0 z \end{bmatrix}$.
Then $\S = [aT']$.
\qed

\begin{lem}\label{noncomm3}
Let $\S \subset \B(\H_1,\H_2)$ be a $1$-hyperreflexive subspace.
Suppose that the orthogonal projection
$Q_{x_0} = P_{\S x_0}$ does not commute with some $Q_y = P_{\S y}$.
Set $Q_0 = (Q_{x_0} \wedge Q_y)^\perp$; and 
let $P_0$ be the projection onto the closed span of all vectors $y$ 
such that $Q_y$ does not commute with $Q_{x_0}$.
Then $P_0x_0 = x_0$ and $Q_0 \S P_0$ is one dimensional, say $\bC T$.

Let $Q$ be the projection onto the range of $T$
and let $P$ be the projection onto the range of $T^*$.
So $T$ is injective on $P\H_1$ with range dense in $Q\H_2$.
The projection $Q$ commutes with $Q_y$ for every $y \in \H_1$;
and $QQ_y$ is $0$, $Q$ or is the rank one projection $[Ty]$. 
Likewise  $P$ commutes with $P_v = P_{\S^* v}$ for every $v \in \H_2$;
and $PP_v$ is $0$, $P$ or is $[T^*v]$. 
\end{lem}

\Prf  By Lemma~\ref{noncomm2}, for any two vectors $y$ and $z$
such that $Q_y$ and $Q_z$ fail to commute with $Q_{x_0}$,
we see that 
\[ Q_{x_0} \wedge Q_y = Q_{x_0} \wedge Q_z = Q_0^\perp \]
and $Q_0 \S |_{\spn\{x_0,y,z\}}$ is one dimensional.  
Moreover for each vector $x$ in $\spn\{x_0,y,z\}$,
$Q_x = Q_0^\perp +[Tx]$.  As $T$ has rank at least two,
it follows that there is a set of vectors $x$ dense in
$\spn\{x_0,y,z\}$ such that $Q_x$ does not commute with $Q_{x_0}$.
Thus the closure of the set of vectors 
$\{x \in \H_1 :   Q_x Q_{x_0} \ne Q_{x_0} Q_x \}$
is a vector space, and thus a subspace.
Also observe that if $Q_y$ does not commute with $Q_{x_0}$,
then $Ty \ne 0$.
Therefore the closure of vectors $x$ so that $Q_x$ does
not commute with $Q_y$ is the same set!

Continuing the analysis of the previous paragraph,
notice that the operator $T$ may be defined on each
subspace $\spn\{x_0,y,z\}$.
The subspace $\bC Tx_0$ does not depend on the choice of
$y$ and $z$.  So it is possible to normalize the choices
by fixing $Tx_0 = u_0$.
Then we obtain a unique value for $Ty$ for each $y = P_0 y$.

Select any element $S \in \S$ such that 
$Q_0^\perp Sx_0  = Tx_0$.
Then $Q_0^\perp Sy$ is a multiple $cTy$ of $Ty$.
In fact it must be exactly $Ty$.
To see this, consider $y_t = (1-t) x_0 + ty$. Then
\[
 Q_0^\perp S y_t = (1-t) Tx_0 + ct Ty.
\]
This has to be a multiple of $Ty_t$ for all $t$.
Since $Tx_0$ and $Ty$ are not collinear,
it follows that $c=1$.

We deduce that $Ty = Q_0^\perp S y$ for all
vectors $y$ such that $Q_y$ does not commute with $Q_{x_0}$.
It follows now that $T$ extends to the closed span
$P_0\H_1$ of these vectors as a bounded operator
with $\|T\| \le \|S\|$.

Now consider a vector $y = P_0 y$ such that $Q_y$ 
commutes with $Q_{x_0}$.
With $y_t$ defined as above, we see that $Ty_t = Q_0^\perp Sy_t \ne
0$ for most values of $t$.
It follows that $Q_{y_t} = Q_0^\perp + [Ty_t]$.
In particular, 
\[
 Q_y \le \limsup Q_{y_t} = Q_0^\perp +[Ty] .
\] 
If $Ty \ne 0$, then $Q_y = Q_0^\perp +[Ty]$.
But if $Ty=0$, we can only deduce that $Q_y \le Q_0^\perp$.

We now define $Q$ to be the projection onto the
range of $T$.  So $Q\le Q_0$ and it commutes with
all $Q_y$ such that $y = P_0 y$.
Similarly, define $P$ to be the projection onto the 
range of $T^*$.

Next suppose that $P_0^\perp z \ne 0$.
Then $Q_z$ commutes with $Q_y$ for
all $y = P_0 y$ for which $Ty \ne 0$.
Therefore it commutes with their intersection
$Q_0^\perp$ and their span $Q_0^\perp + Q$.
Therefore $Q_z Q$ is a  projection which
commutes with $[Ty]$ for all $y = P_0 y$.
As the range of $T$ is dense in $Q\H_2$,
it follows that $Q_z Q$ is either 0 or $Q$.

Similarly, consideration of $\S^*$ shows that if
there are vectors $u$ and $v$ such that the range
projections $P_u$ and $P_v$ onto $\ol{\S^*u}$ and $\ol{\S^*v}$
do not commute, one likewise finds projections
$P_0$ and $Q_0$ so that $P_0 \S^* Q_0 =\bC P_0 T^* Q_0$.
However once one finds such a form, one also sees
that the projections $Q_x$ and $Q_y$ for $x,y \in P_0\H_1$
would also fail to commute, and that we have already identified
this subspace in the previous analysis.

What we can conclude is that for $u = P_0 u$
such that $T^*u \ne 0$, $P_u = P_0^\perp + [T^* u]$;
and if $T^*u=0$, then $P_u \le P_0^\perp$.
Also if $P_0^\perp v \ne 0$, then $P_v$ commutes
with $P_0^\perp$ and $P$, and $P_v P$ is either 0 or $P$.
\qed

\begin{lem}\label{noncomm4}
Let $\{(P_i,Q_i) : i \in \I \}$ be the collection of all pairs of 
projections $P \in \B(\H_1)$ and $Q\in \B(\H_2)$ obtained as in
Lemma~$\ref{noncomm3}$ from a pair of vectors $x,y$ such that
$Q_x$ and $Q_y$ do not commute.  
Then $\P = \{ P_i : i \in \I \}$ and  $\Q = \{ Q_i : i \in \I \}$
are families of pairwise orthogonal projections which commute with
every $P_x$ and $Q_x$ respectively.
For $x \in \H_1$, there is at most one $i \in \I$ such that 
$Q_i Q_x$ is neither $0$ nor $Q_i$.
Moreover, $\S$ is a $\fC$--$\fD$ bimodule, where $\fC$ and $\fD$ are the
abelian von Neumann algebras generated by the $\{P_i\}$ and
$\{Q_i\}$ respectively.
\end{lem}

\Prf As before, we write $Q_x = [\S x]$ and $P_u = [\S^* u]$.
For each pair of vectors $x,y \in \H_1$ such that $Q_x$ and $Q_y$
do not commute, Lemma~\ref{noncomm3} provides projections 
$P \in \B(\H_1)$ and $Q\in \B(\H_2)$ so that $Q\S P = \bC T$
is 1-dimensional. Moreover every $Q_z$ commutes with $Q$; and
$Q_z Q$ is 0, $Q$ or $[Tz]$. 

This immediately implies that if $x',y'$ is another such pair,
then the corresponding projections $P'$ and $Q'$ either equal 
$P$ and $Q$ or they are orthogonal.  
Thus there is a set
$\{(P_i, Q_i) : i \in \I \}$ consisting of all such pairs.
We write $Q_i \S P_i = \fX_i = \bC T_i$.

Let $x \in \H_1$.  
Suppose that for some $i\in \I$, $Q_iQ_x = [T_i x] \ne 0$.
Then there is another vector $y$ so that $Q_i Q_y = [T_i y]$
does not commute with $[T_i x]$.  Hence by Lemma~\ref{noncomm1},
$Q_x = (Q_x \wedge Q_y) + [T_i x]$.
Moreover, if $z$ is any other vector such that $Q_z$ does not
commute with $Q_x$, then $Q_z = (Q_x \wedge Q_y) + [T_i z]$.
Hence for all $z = P_j z$ where $j \ne i$, $Q_z$ commutes with
$Q_x$ and indeed with all $Q_y$ for which $Q_i Q_y = [T_i y] \ne
0$.  It follows that $Q_z Q_i$ is 0 or $Q_i$, not a one dimensional
projection.  Likewise $Q_x Q_j$ is 0 or $Q_j$.

In particular, $Q_i$ commutes with every $Q_x$.
It follows that if $S \in \S$ and $x \in \H_1$, then
\[
 Q_i S x \in Q_i Q_x \H_2 \subset Q_x \H_2 = \ol{\S x} .
\]
By the reflexivity of $\S$, $Q_i S \in \S$.
Consequently $\fC \S = \S$.
Consideration of $\S^*$ yields a similar conclusion on the right.
\qed

\begin{lem}\label{reflexiveX}
Let $\S$ be a 1-hyper\-reflexive subspace of $\B(\H_1,\H_2)$.
Let $\{(P_i,Q_i) : i \in \I \}$ be the collection of all pairs of 
projections obtained in Lemma~$\ref{noncomm4}$.
Define a space of operators by
\[
 \fX = \S + \wotsum_{i \in \I} Q_i \B(\H_1,\H_2) P_i .
\]
Then $\fX$ is reflexive and the projections $R_x := [\fX x]$
commute with $R_y$ and $Q_y$ for $y \in \H_1$.
\end{lem}

\Prf
Suppose that $X\in\Ref(\fX)$ and let $D=\sum_{i\in \I} Q_iXP_i$.
Then $D\in \fX$, and to show that $X\in \fX$ it suffices to prove
that $X-D\in \S$.
Suppose that $x \in \H_1$ is a vector such that $Q_i Q_x \in \{ 0,
Q_i\}$ for all $i \in \I$.
Then in particular, if $P_i x \ne 0$, since $T_i$ is
injective on $P_i\H_1$, $Q_i Q_x \H_2$ contains the non-zero
vector $T_i P_i x$; and so $Q_x Q_i = Q_i$.
Therefore
\[
 \fX x = \S x + \sum_{\{i: P_i x \ne 0 \}}Q_i\H_2 = \S x .
\]
Otherwise there is a unique $i_0 \in I$ for which 
$Q_{i_0}Q_x = [T_{i_0}P_{i_0}x] \ne 0$.
As before, for all other $i$ for which $P_i x \ne 0$, 
$Q_x Q_i = Q_i$.
Also $Q_{i_0} \S P_{i_0}^\perp x = 0$ (since otherwise it is
$Q_{i_0} \H_2$). 
Thus
\[
  \fX x = \S x + \sum_{\{i: P_i x \ne 0 \}}Q_i\H_2 
        = \S x + Q_{i_0}\H_2 .
\]
Since $\S$ is a $\fC$--$\fD$ bimodule, 
\[
 \S P_{i_0}^\perp x = Q_{i_0}^\perp \S P_{i_0}^\perp x
  \subset Q_{i_0}^\perp \S x \subset \S x.
\]
Therefore
\[
 (X-D) x = (X-D) P_{i_0}^\perp x + (X-D) P_{i_0} x
         \in Q_{i_0}^\perp \S x \subset \S x .
\]
Hence $X-D$ belongs to $\S$ as claimed.

Since $Q_{i_0}$ commutes with $Q_x$, 
\[
 R_x = Q_x \vee Q_{i_0} = Q_x + Q_{i_0} - Q_{i_0}Q_x .
\]
$R_x$ therefore commutes with all $Q_j$.
If $y \in \H_1$, either $Q_y$ commutes with $Q_x$ and so with
$R_x$ or $Q_y = Q_x - [T_{i_0}x] + [T_{i_0}y]$ which evidently also
commutes with $R_x$.
Finally we can conclude that $R_x$ also commutes with $R_y$.
\qed

\subsection*{Proof of Theorem~\ref{T:1hyp}.} \strut\\
Let $\S$ be a 1-hyper\-reflexive subspace of $\B(\H_1,\H_2)$.
As in Lemma~\ref{noncomm4},
define $\P = \{ P_i : i \in \I \}$ and  $\Q = \{ Q_i : i \in \I \}$,
and select operators $T_i$ such that $Q_i\S P_i=\bC T_i$.  
Set
\[
 \fX = \S + \wotsum_{i \in \I} Q_i \B(\H_1,\H_2) P_i 
\]
as in the previous lemma.

Let $\fC$ be the abelian von Neumann algebra generated by the set
$\Q \cup \{ R_x : x \in \H_1 \}$, and let $\fD$ be generated by 
$\P \cup \{ [\fX^* u] : u \in \H_2 \}$. 
As in the proof of Theorem~\ref{masa modules},
the reflexive subspace $\fX$ is a $\fC'$--$\fD'$ bimodule.
 From the proof of Lemma~\ref{reflexiveX}, we see that $\Q$ are
atoms of $\fC$ and $\P$ are atoms of $\fD$.

Following the proof of Theorem~\ref{masa modules}, we obtain 
families of projections $\C = \{ C_j : j \in \J \}$ and  
$\D = \{ D_j : j \in \J \}$ so that  
\[
 \fX = \{ T \in \B(\H) : C_jTD_j \in \fX_j \qforal j \in \J \} 
\]
where $\fX_j$ are $\fC$--$\fD$ bimodules in $C_j \B(\H_1,\H_2) D_j$
for $j \in \J$.

Furthermore, as in that proof, each $\fX_j$ is 1-hyper\-reflexive if and
only if it is a nest bimodule.  
Corollary~\ref{C:masa} shows that failure to be 1-hyper\-reflexive yields
orthogonal projections $D_1, D_2, D_3 \in \fD$ and 
$C_1, C_2 \in \fC$ so that $(C_1+C_2) \fX (D_1+D_2+D_3)$ has the form
\[
 \begin{bmatrix} * & 0 & 0\\0 & * & 0\end{bmatrix}
\]
where the $*$ entries are non-zero; or the analogous $3 \times 2$ form.

We claim that there are $i_1, i_2 \in \I$ so that $P_{i_k} \le D_k$
and $Q_{i_k} \le C_k$ for $k = 1,2$.
If not, then 
\[ (C_1+C_2) \S (D_1+D_2+D_3) = (C_1+C_2) \fX (D_1+D_2+D_3) \]
is not 1-hyper\-reflexive, contrary to hypothesis.
By cutting down, we may suppose that $D_k = P_{i_k}$
and $C_k = Q_{i_k}$ for $k = 1,2$.
Thus
\[
 (C_1+C_2) \S (D_1+D_2+D_3) = 
 \begin{bmatrix} aT_{i_1} & 0 & 0\\0 & bT_{i_2} & 0\end{bmatrix} .
\]

By Theorem~\ref{compress}, this is 1-hyper\-reflexive.
And so Theorem~\ref{T:22case} shows that $a$ and $b$ are dependent,
say $b = ar$.
But then this compression consists of multiples of the operator
$\begin{bmatrix}T & \begin{matrix}0\\0 \end{matrix}\end{bmatrix}$ where 
$T = \begin{bmatrix}T_{i_1} & 0 \\0 & T_{i_2}\end{bmatrix}$.
However in this case, any vector $x = (D_1+D_2)x$ satisfies
$(C_1+C_2)Qx = [Tx]$.
Recall that $T_{i_k}$ are injective on $D_i \H_1$ and have range dense in
$C_k\H_2$.  Hence $[Tx]$ is non-zero for all non-zero $x = (D_1+D_2)x$.
It is easy to see that if $x_k = D_k x_k$ are non-zero, then
$Q_{x_1}$ and $Q_{x_1+x_2}$ cannot commute.
This contradicts the construction of the projections $D_i$ by
Lemma~\ref{noncomm4}.
Consequently, $a$ and $b$ are independent and the compression is not 
 1-hyper\-reflexive.
As this is not possible, we deduce that $\fX_j$ is a nest bimodule.

Finally, observe that $C_j \S D_j$ is obtained from $\fX_j$ by the
construction of Lemma~\ref{tri_const}.
We have seen that $\fX_j$ is a nest bimodule and that $Q_i \S P_i =
\bC T_i$.
Let $\N$, $\M$ be nests such that $\fX_j$ is a $\T(\N)$--$\T(\M)$
bimodule, and let $\theta:\N\to\M$, $\theta^*:\M\to\N$ be the
functions so that $\fX_j=\fX(\theta)$ and $\fX_j^*=\fX(\theta^*)$. We
may assume that $\N$ and $\M$ are minimal in the sense that
$\N=\theta^*(\M)$ and $\M=\theta(\N)$.
Define nests $\N'$ and $\M'$ by
\begin{align*}
\N'&=\N\cup\{N+P_i: N\in \N,\ i\in \I,\ N<N+P_i<N^+\},\\
\M'&=\M\cup\{M^+-Q_i: M\in \M,\ i\in \I,\ M<M^+-Q_i<M^+\} 
\end{align*}
and let $\fX_j^0=\{X\in \fX_j: Q_iXP_i=0\FOR i\in \I\}$. Then
$(\fX_j,\fX_j^0)$ is a $\T(\N')$--$\T(\M')$ bimodule pair, and 
$\Delta(\fX_j,\fX_j^0)=\{(P_i,Q_i):i\in \I\}$.
It remains to verify that the subspaces $Q_i\S P_i=\bC T_i$ are independent.
To that end, it would suffice to show that $Q_i \S P_i$ belongs to $\S$
for each $i \in \I$. 
This follows from the reflexivity of $\S$.
For if $P_i x = 0$, $T_i x = 0 \in \S x$.
Suppose that $x_i :=P_ix \ne 0$, say $x = x_i + y$.
Then $Q_i Q_y$ is either $0$ or $Q_i$ while $Q_i Q_{x_i} = [T_i x_i]$
is exactly one dimensional.
If $Q_i Q_y = 0$, then $Q_i Q_x = [T_i x_i]$; while if 
$Q_i Q_y = Q_i$, then $Q_i Q_x = Q_i$.
In either case, $\S x$ contains $T_ix = T_i x_i$.
By the reflexivity of $\S$, $T_i \in \S$ as desired.
\qed

\section{1-Hyperreflexive Algebras}\label{S:1hypalg}

In this section, we apply Theorem~\ref{T:1hyp} to the 
case of unital algebras.

\begin{cor}\label{C:1hypalg}
A unital algebra~$\fA$ is a $1$-hyper\-reflexive if and only if either
\begin{enumerate}
\item there is a nest $\N$ and a collection $\{ A_i : i \in \I \}$ 
of atoms of $\N$ such that
\[ \fA = \{ T \in \T(\N) : A_i T A_i \in \bC A_i \FORAL i \in \I \} \]
or
\item there is a projection $P$ so that $\fA = \B(P\H) \oplus \B(P^\perp\H)$.
\end{enumerate}
\end{cor}

\Prf The examples of (1) are unital algebras which are
1-hyper\-reflex\-ive as subspaces by Lemma~\ref{tri_const},
and hence as algebras; while (2) is example A3.
Both fall under the rubric of Theorem~\ref{T:1hyp}.

Conversely, consider the construction of Lemma~\ref{diag_const}
and suppose that there are two or more diagonal blocks involved.
Then there are projections $C_1,C_2$, $C_3 = (C_1+C_2)^\perp$,
$D_1,D_2$ and $D_3 = (D_1+D_2)^\perp$ so that 
$\fA$ contains 
\[
 \{ T \in \B(\H) : C_iTD_i = 0,\, i = 1,2,3 \}
\]
where, since $C_3$ and $D_3$ would consume all but the first two
blocks, this includes all cases in which there are at least two such
blocks if we allow one or both of $C_3$ and $D_3$ to be 0.
So 
\[
 \fA = \{ T \in \B(\H) : C_iTD_i \in \fX_i,\, i = 1,2,3 \}
\]
where $\fX_i = C_i \fA D_i$.  
Moreover, $\fX_1$ and $\fX_2$ are proper subspaces of the form of
Lemma~\ref{tri_const}.
$\fX_3$ is an arbitrary 1-hyperreflexive subspace.

Suppose that $(D_2+D_3)(C_2+C_3)\ne 0$.
Then 
\begin{align*}
 \fA D_1 &\supset \fA D_1 + \fA (D_2+D_3) \fA D_1 \\
 &\supset  (C_2+C_3) \B(\H) D_1 + \fA (D_2+D_3) (C_2+C_3) \B(\H) D_1
 \\ &\supset (C_2+C_3) \B(\H) D_1 + C_1 \B(\H)(D_2+D_3) (C_2+C_3) \B(\H) D_1
\end{align*}
The \wot-closed span of $\B(\H)(D_2+D_3) (C_2+C_3) \B(\H)$ is
all of $\B(\H)$.  
As $\fA D_1$ is a \wot-closed subspace, we would conclude that
$\fA D_1 = \B(\H) D_1$.
However $C_1\fA D_1 = \fX_1$ is proper, so $D_2 + D_3 \le C_1$.
Similarly $D_1 + D_3\le C_2$.
Consequently $D_3=0$.
Therefore $D_2 = C_1$ and $D_1 = C_2$, and $C_3=0$.

If $\fX_1=\fX_2=0$, this yields $\fA = \B(C_1\H) \oplus \B(C_2\H)$
which is case (2).  Otherwise, since $\fA$ is a \wot-closed algebra
containing this as a subalgebra, $\fX_1 = \B(C_2\H,C_1\H)$ or 0, and
$\fX_2 = \B(C_1\H,C_2\H)$ or 0.  
The three possibilities are all nest algebras.

Also, there is the case in which there is one block 
$\fX_1 \subset C_1\B(\H) D_1$, but at least one of $C_1,D_1$ is a
proper projection.  If both $C_2$ and $D_2$ are non-zero, the same
argument shows that $\fA = \B(\H)$.  If $C_2=0$ and $D_2 \ne 0$,
then we may suppose that $D_2$ is the largest subspace such that $\fA
D_2 = \B(\H) D_2$. Then observe that if $D_2 \fA D_1 \ne 0$, there
are vectors $x = D_1 x$ such that $\fA x = \H$, contrary to our
assumption on $D_2$.  Thus $\fA D_1 = D_1\H$.  As $\fX_1 = \fA D_1$
has the form of Lemma~\ref{tri_const}, it is clear that 
$\fA = \fX_1 + \B(\H)D_2$ also has this form.
A similar analysis holds if $D_2=0$.
Hence we may now assume that the 1-hyper\-reflexive unital
algebra $\fA$ consists of a single block obtained using
Lemma~\ref{tri_const}.

We have reduced the problem to the situation where there are
nests $\M$ and $\N$ and a bimodule pair $(\fX, \fX_0)$ with atoms 
$\{(A_i,B_i) : i \in \I \}$ and operators $X_i \in \B(B_i\H,A_i\H)$
so that 
\[
 \fA = \{ T \in \fX : A_iTB_i \in \bC X_i,\, i \in \I \} 
\]
and
\[
 \fX_0 = \{ T \in \fX : A_iTB_i = 0,\, i \in \I \} .
\]
Let $M_i$ and $N_i$ be the elements of $\M$ and $\N$
respectively such that $A_i = M_i^+-M_i$ and $B_i = N_i^+-N_i$ for
$i \in \I$.

Let $\theta$ and $\theta_0$ be the left continuous order
preserving maps of $\N$ into $\M$ such that $\fX = \fX(\theta)$
and $\fX_0 = \fX(\theta_0)$.
Recall that $\theta_0(N) = \theta(N)$ unless $N = N_i^+$ for
some $i \in \I$, in which case 
$\theta_0(N_i^+) = M_i \supset \theta(N_i)$ and $\theta(N_i^+) = M_i^+$. 
We may assume that each $X_i$ is injective on $B_i \H$ with
range dense in $A_i\H$.  In particular, this ensures that  
$\ol{\fA N} = \theta(N)$ for all $N \in \N$.

Observe that if $N_i^+ \subset M_i$, then
\[
 M_i \supset \theta(N_i) = \ol{\fA N_i} = \ol{\fA^2 N_i} = \ol{\fA M_i} 
 \supset \ol{\fA N_i^+} = M_i^+ ,
\]
which is absurd.  If $x = B_i x \not\in M_i$, then because $\fA$
is unital,
\[
 \ol{\fA x} = M_i \vee \bC X_i x \supset M_i \vee \bC x .
\]
Thus $\ol{\fA x} = M_i \vee \bC x$.  By lower semicontinuity, this
identity persists for all $x = B_i x$.
Since $X_i$ is injective, this means that $M_i \cap B_i\H = \{0\}$.
Now 
\[
 M_i \vee A_i\H = M_i^+ = \ol{\fA N_i^+}
 = \bigvee_{x = B_i x} \fA x = M_i \vee B_i\H .
\]
That is, $M_i^+ = M_i \vee B_i\H$.

Suppose that $M \in \theta(\N)$, and let $N$ be the largest element of
$\N$ with $\theta(N) = M$, which exists since $\theta$ is left
continuous.  Also let $N'$ denote the smallest element of $\N$
containing $M$. 
Then since $\fA$ is unital, $N \subset \ol{\fA N} = M \subset N'$.
So if $N'=N$, we see that $M=N$.  Otherwise
\[
  M =  \ol{\fA N} = \ol{\fA^2 N} = \ol{\fA M}
  \supset \ol{\fX_0 M} = \fX_0 N' = \theta_0(N') .
\]
Since $N' > N$, $M':=\theta(N') > M$ and so 
\[
 M = \theta_0(N') \ge \theta_0(N')^- = (M')^- .
\]
Therefore $M = (M')^-$.
It follows that $N' = N_i^+$ for some $i \in \I$;
and hence $N=N_i$, $M=M_i$ and $M' = M_i^+$.
But $N_i \subset M_i$, and by the previous paragraph, $M_i \cap B_i\H =
\{0\}$.
Since $M_i \subset N_i^+$, this means that $M_i=N_i$ and in fact
$N'=N_i$ after all.

Consequently we deduce that whenever $N$ is the largest element of
$\N$ with $\theta(N) = M \in \M$, then $M=N$.
This includes every $N_i$, $i \in \I$ because
$\theta(N_i^+)>\theta(N_i)$. 
Again the analysis of the atoms yields
\[
 \theta(N_i^+) = M_i \vee B_i\H = N_i \vee B_i\H = N_i^+.
\]
Since $\theta(N) \ge N$, we see that each $N_i^+$ is also the largest
element of $\N$ with $\theta(N) = N_i^+$.

Let $\N_0 = \{N \in \N : \theta(N) = N \}$.
We claim that $\N_0$ is complete.  
Indeed, if not, there is an element $N$ in the completion which is a
\textit{monotone} limit of elements $N_\alpha \in \N_0$.
If it is a limit from below, then by left continuity
\[
 \theta(N) = \sup \theta(N_\alpha) = \sup N_\alpha = N .
\]
And if it is a limit from the above,
\[
 N \le \theta(N) \le \inf \theta(N_\alpha) = \inf N_\alpha = N .
\]
Moreover $\N_0$ contains $\theta(\N)$; whence $\N_0 = \theta(\N)$.

Now we have $\fX = \T(\N_0)$ and $\fX_0$ is the ideal
\[ \{T \in \T(\N_0) : A_iTA_i=0,\, i\in\I \} .\]
Since $\fA$ is unital, it follows that each $X_i = A_i$ for $i \in \I$.
So $\fA$ has the form of (1).
\qed

\section{Complete Hyperreflexivity}

As mentioned in the introduction, we make the following natural
definition.  
Here for $\S \subset \B(\H)$, we denote by $\S \botimes \B(\K)$ the
\wot-closure of the spatial tensor product in $\B(\H \otimes \K)$.

\begin{defn}\label{D:comp_hyp}
A \wot-closed subalgebra (or subspace) $\S$ of $\B(\H)$ is 
\textit{completely hyperreflexive} if $\S \botimes \B(\K)$
is hyperreflexive for a separable Hilbert space $\K$.
Let $\kappa_\S^c$ denote the hyperreflexivity constant of 
$\S \botimes \B(\K)$.
\end{defn}

Unlike an arbitary compression considered in Theorem~\ref{compress}, it
is a very different situation when the compression remains in the
subspace. This lemma also appears in~\cite{KL} with a different
proof.

\begin{lem}\label{L:cutdown}
Let $\S$ be a subspace of $\B(\H_1,\H_2)$.
Suppose that $P$ and $Q$ are projections such that $Q \S P \subset \S$.
Then considering $Q \S P$ as a subspace of $\B(P\H_1, Q\H_2)$, we
obtain the inequality $\kappa_{Q \S P} \le \kappa_\S$. 
\end{lem}

\Prf Let $T\in \B(P\H_1, Q\H_2)$. Then $\dist(T,Q\S P)=\dist(T,\S)$
and 
\begin{align*}
  \beta_{\S}(T) &= \sup_{\|x\|=1}\inf_{S\in\S} \|Tx-Sx\|
  \leq \sup_{\|x\|=1}\inf_{S\in Q\S P}\|Tx-Sx\|\\
  &=\sup_{\substack{\|x\|=1,\\ x=Px}} \inf_{S\in Q\S P}\|Tx-Sx\|
  =\beta_{Q\S P}(T).
\end{align*}
So 
\begin{align*}
\kappa_{Q\S P} &= \sup \{ 1/\beta_{Q\S P}(T) : T=QTP,\ \dist(T,Q\S P)=1\} 
\\ &\leq \sup \{ 1/\beta_{\S}(T) : T=QTP,\ \dist(T,\S)=1\} \leq
\kappa_{\S}.
\end{align*}
\upqed

\begin{prop}\label{P:increasing}
If $\S$ is a \wot-closed subspace, then the hyperreflexivity constants for
$\S \otimes \fM_n$ are increasing, and
\[ \kappa_\S^c := \dlim_{n\to\infty} \kappa_{\S \otimes \fM_n} .\]
\end{prop}

\Prf Fix an orthonormal basis $\{e_n\}_{n\ge1}$ for $\K$. For $n\ge
1$, let $P_n$ be the orthogonal projection $I_{\H}\otimes Q_n\in
\B(\H\otimes K)$ where $Q_n\in \B(\K)$ is the orthogonal projection
onto $\spn\{e_1,e_2,\dots,e_n\}$. We identify $\S\otimes \fM_n$
and $P_n(\S\botimes \B(\K))$.  Lemma~\ref{L:cutdown} applies to show
that
\[
 \kappa_{\S \otimes \fM_n} \le \kappa_{\S \otimes \fM_{n+1}}
 \le \kappa_\S^c .
\]
Denote the limit by $\kappa = \dlim_{n\to\infty} \kappa_{\S \otimes\fM_n}$. 
Fix $T \in \B(\H\otimes\K)$ and for $n\ge1$, let $T_n=P_nTP_n$.
Then
\begin{align*}
 \beta_{\S \botimes \B(\K)}(T) &\ge
 \sup_{\|x\| = \|P_n x\| = 1}\inf_{S \in \S \botimes \B(\K)} \| (T - S)x \|
 \\&\ge \sup_{\|x\| = \|P_n x\| = 1}
  \inf_{S \in \S \botimes \B(\K)} \| P_n (T- S) P_n x \|
  \\&= \sup_{\|x\| = \|P_n x\| = 1}
  \inf_{S \in \S \otimes \fM_n} \| (T_n - S) x \|
 \\&\ge \kappa^{-1} \dist(T_n, \S \otimes \fM_n) .
\end{align*}
It is easy to verify from the lower semicontinuity of the norm in
the strong operator topology that 
\[
  \lim_{n\to\infty} \dist(T_n, \S \otimes \fM_n) = \dist(T, \S \botimes  \B(\K)).
\]
Thus taking a supremum over all $n\ge1$ in the previous expression yields
\[
 \beta_{\S \botimes \B(\K)}(T) \ge \kappa^{-1} \dist(T,\S \botimes  \B(\K)) .
\]
So $\kappa^c_\S \le \kappa \le \kappa^c_\S$; whence equality holds.
\qed

A long-standing question posed in \cite{DO} is:

\begin{ques} 
Is every hyperreflexive subspace completely hyperreflexive?
\end{ques}

\begin{prop}\label{two}
The algebra $\D \subset \fM_3$ consisting
of matrices of the form 
$\begin{bmatrix}\alpha&0&0\\0&\beta&0\\0&0&\beta\end{bmatrix}$
has distance constant at least $2/\sqrt3$.
\end{prop}

\Prf Consider the matrix
\[
 T = 
 \begin{bmatrix} 
  \phantom{-}0&\phantom{-}0&\sqrt2\,\\
  -\sqrt2&-1&0\\
  \phantom{-}0&\phantom{-}0&1
 \end{bmatrix} .
\]
Observe that 
\begin{align*}
 \left\| \begin{bmatrix} 
  \phantom{-}\alpha&0&\sqrt2\\
  -\sqrt2&\beta-1&0\\
  \phantom{-}0&0&\beta+1
 \end{bmatrix} \right\| &\ge
 \max \Big\{ \| \begin{bmatrix}-\sqrt2&\beta-1\end{bmatrix} \|,\
  \left\| \begin{bmatrix} \sqrt2\\ \beta+1 \end{bmatrix} \right\|
  \Big\} \\
  &= \max \{ \big( 2+ |\beta \pm 1 |^2 \big)^{1/2} \}.
\end{align*}
The minimum over all $\alpha,\beta$ is $\sqrt3$, and this is attained
when $\alpha = \beta = 0$.

Next we note that the proper invariant subspaces have the form
\[
 \spn\{e_1\},\ \spn\{e_1\}^\perp, \spn\{v_s\}, \qand \spn\{v_s\}^\perp
\]
where $v_s = (0,c,s)^t$, $|s| \le 1$, and $c = \sqrt{1-|s|^2}$.
If $V$ is a 1-dimensional invariant subspace containing a unit vector
$v$, we can compute 
\[
 \| P_V^\perp T P_V \| = \| P_V^\perp T v \| = \| Tv - \ip{Tv,v}v \| .
\]
While if $P$ is a 2-dimensional invariant subspace orthogonal to a unit
vector $v$, we can instead compute
\[
 \| P_V^\perp T P_V \| = \| T^*v - \ip{T^*v,v}v \| .
\]
For $\spn\{e_1\}$ and $\spn\{e_1\}^\perp$, we obtain $\sqrt2$.

Consider $V = \spn\{v_s\}$.  Then
\begin{align*}
 \| P_V^\perp T P_V \| &=
 \left\| \begin{bmatrix}\sqrt2 s\\-c\\s\end{bmatrix} 
   - \big( |s|^2-c^2 \big) \begin{bmatrix}0\\c\\s\end{bmatrix} \right\| 
  = \left\| \begin{bmatrix}\sqrt2 s\\-2s^2c\\2sc^2\end{bmatrix} \right\|
 \\ &= \sqrt{ 2|s|^2 + 4 |s|^2 (1-|s|^2) } \le \frac 3 2 .
\end{align*}
This bound is attained when $s = \sqrt3/2$.
Thus $\beta_\D(T) = 3/2$.

Therefore $\kappa_\D \ge \frac{\sqrt3}{3/2} = \frac 2 {\sqrt3}$.
\qed

\begin{cor}\label{dim2}
If $\S$ is a subspace of $\fM_n$ or $\B(\H)$ of dimension at least
$2$, then $\S \otimes \bC I_n$ is not $1$-hyperreflexive for
any $n\ge3$ or $n=\infty$.
\end{cor}

\Prf Choose unit vectors $x_1, x_2$ and $y_1,y_2$ so that
the functionals $\psi_i(A) = \ip{A x_i,y_i}$ on $\S$ are linearly
independent.
In $\H \otimes l^2_n$, let $P$ be the projection onto 
$\spn\{ x_1 \otimes e_1, x_2 \otimes e_2, x_2 \otimes e_3 \}$ and let $Q$ be the projection onto
$\spn\{  y_1 \otimes e_1, y_2 \otimes e_2, y_2 \otimes e_3 \}$.
Consider the compression $Q\S P$.  
This is evidently the algebra $\D$ of Proposition~\ref{two}.
Since $\D$ does not have distance constant 1, 
Theorem~\ref{compress} shows that  
$\S \otimes \bC I_n$ also does not have distance constant 1.
\qed

\begin{cor}\label{scalars}
The algebra $\bC I$ of scalar matrices is completely hyperreflexive.
However $1 = \kappa_{\bC I} < \kappa^c_{\bC I}$.
\end{cor}

Getting an explicit lower bound greater than 1 takes a lot more work.
As far as we know, the bound that we can get is not very good.

\begin{prop}
  The complete distance constant $\kappa$ for the algebra $\bC I$ is at
  least $1.03$.
\end{prop}

\Prf Let $\S=\bC I\otimes \fM_{2,1}\subset \B(\H)\otimes \fM_{2,1}$.
By Lemma~\ref{L:cutdown}, $\kappa \ge \kappa_\S$.

Let $\alpha=\sin(\pi/8)$ and $\beta=\cos(\pi/8)$. 
Let $\K$ be a two-dimensional subspace of $\H$
with orthonormal basis $\{e_1,e_2\}$.
Define $T_1,T_2\in \B(K)$ by
\[
 T_1 =
 \begin{bmatrix}1/\sqrt2\\ 1/\sqrt2\end{bmatrix}
 \begin{bmatrix}\alpha & -\beta \end{bmatrix} =
 \frac 1 {\sqrt2} \begin{bmatrix}
    \alpha & -\beta \\ \alpha & -\beta
  \end{bmatrix}
\]
and
\[
 T_2 =
 \begin{bmatrix}0\\ 1 \end{bmatrix}
 \begin{bmatrix}\beta & \alpha \end{bmatrix} =
 \begin{bmatrix}
    0 & 0 \\ \beta & \alpha
 \end{bmatrix} .
\]
Consider $\H \oplus \H$ as $\H \otimes \bC^2$ and let $\{u_1,u_2\}$
be the standard basis for $\bC^2$.
Decompose $\H = \K\oplus \K^\perp$ and
\[
 \H \oplus \H = (\K\otimes u_1) \oplus (\K^\perp\otimes u_1) \oplus
                (\K\otimes u_2) \oplus (\K^\perp\otimes u_2)
\] 
Define $T\in \B(\H)\otimes \fM_{2,1} = \B(\H,\H\oplus \H)$ by
\[
T=  \begin{bmatrix}
    T_1 & 0 \\ 0 & 0 \\ T_2 & 0 \\ 0 & 0
  \end{bmatrix}.
\]
Let $f_1 = \tfrac{1}{\sqrt2}(e_1+e_2)\otimes u_1$ and $f_2 =
e_2 \otimes u_2$.
Let $P$ be the orthogonal projection onto $\K$ and let
$Q$ be the orthogonal projection onto $\spn\{ f_1,f_2 \}$,
which is the range of $T$. 
With respect to the bases $\{e_1,e_2\}$ and $\{f_1,f_2\}$,
\[
  QTP = 
  \begin{bmatrix}
    \alpha & -\beta \\ \beta & \alpha
  \end{bmatrix}
\qand
 Q\S P = \left\{
  \begin{bmatrix}
    a & a \\ 0 & b 
  \end{bmatrix}
 : a,b\in \bC \right\}.
\] 
Observe that $QTP$ is unitary.
The proof of Theorem~\ref{T:22case} shows that
$\dist(QTP, Q\S P) = 1 = \|T\|$. As in the proof of
Theorem~\ref{compress}, it follows that $\dist(T,\S) = 1$.
Thus $\beta_\S(T)^{-1} \le \kappa_\S \le \kappa$.

Let $v = xe_1 + ye_2 + z$ be a unit vector in $\H$ where $x,y\in\bC$ and
$z \in \K^\perp$. 
Then $\S v = v \otimes \bC^2$, and a computation shows that
\begin{align*}
  \| P_{\S v}^\perp Tv \|^2  &= 
  \| Tv \|^2 - |\ip{Tv, v\otimes u_1}|^2 - |\ip{Tv, v\otimes u_2}|^2 \\
  &= |x|^2 + |y|^2 - \tfrac12 | x + y |^2 | \alpha x - \beta y|^2 
   - |y|^2 |\beta x+\alpha y|^2 .
\end{align*}
If we expand this expression and consider what happens when $|x|$,
$|y|$ and $\|z\|$ are fixed, then the only variable terms form a
quadratic in $\rho = 2\re(x\ol{y})$ whose leading term is
$\tfrac12\alpha\beta \rho^2$. Since $\alpha\beta>0$, this function is
maximized over the interval $\rho\in [-2|xy|, 2|xy|]$ at an endpoint,
so we may assume that $x$ and $y$ are real.

Let 
\[
 \psi(x,y) =
 \tfrac12 (x+y)^2 (\alpha x - \beta y)^2 + y^2 (\beta x+ \alpha y)^2,
\]
so that the expression above becomes
\[
 \|P_{\S v}^\perp Tv\|^2 = x^2 + y^2 - \psi(x,y).
\]
Let $k= \displaystyle\inf_{x^2+y^2=1}\psi(x,y).$ 
Numerical experiments reveal that $k > 0.058$.
Since $\psi(rx,ry) = r^4\psi(x,y)$ and $k < 1$,
\begin{align*}
  \beta_{\S}(T)^2 & =\sup_{\substack{0\leq r\leq 1\\ x^2+y^2=r^2}}  \|P_{\S v}^\perp
    Tv\|^2\\
    &= \sup_{0\leq r \leq 1} \big( r^2 -
       \inf_{x^2+y^2=r^2}\psi(x,y)\big)\\
    &= \sup_{0\leq r \leq 1} r^2 - kr^4 = 1 - k < 0.942.
\end{align*}
Thus $\kappa \ge (1-k)^{-1/2} > 1.03$. 
\qed

Now $\bC I \otimes \B(\H)$ is a type I von Neumann algebra.
So it is hyperreflexive by Christensen's Theorem \cite{Chr1}
(see \cite[Theorem~9.6]{D1}) with constant at most 4.
We will show that the constant is at most 2.

\begin{prop}\label{P:scalartensor}
The distance constant for $\bC I \otimes \B(\H)$ is at most $2$.
The distance constant for $\bC I_2 \otimes \B(\H)$ is at most
$\frac32$.
\end{prop}

\Prf The idea is to average over a group of unitaries which have
two point spectrum.
Consider the eight element group $\G$ consisting of matrices
\[
 \begin{bmatrix}\pm1&0\\0&\pm1\end{bmatrix} \qand
 \begin{bmatrix}0&\pm1\\ \pm1&0\end{bmatrix} .
\]
Observe that $\G$ has trivial commutant.
Hence $\G \otimes I$ has commutant $\bC I_2 \otimes \B(\H)$.
Also note that the spectrum of an element of $\G \otimes I$ is one of
$\{1\}$, $\{-1\}$, $\{\pm1\}$ or $\{\pm i\}$.
So the elements can be written as $\pm I_2$, $2P-I_2$ or
$i(2Q-I_2)$ where $P \otimes I$ and $Q \otimes I$ are projections in 
$\fM_2 \otimes \bC I$, and hence in $\Lat(\bC I_2 \otimes \B(\H))$.
Indeed $P$ is one of $E_{11}\otimes I$, $E_{22}\otimes I$ or 
$\begin{bmatrix}\frac12 I&\pm\frac12 I\\
    \pm\frac12 I&\frac12 I\end{bmatrix}$
and $Q$ is one of $\begin{bmatrix}\frac12 I&\pm\frac i2 I\\
    \mp\frac i2 I&\frac12 I\end{bmatrix}$.

Define an expectation onto $\bC I_2 \otimes \B(\H)$ by
\begin{align*}
  \Phi(T) &= \frac18 \sum_{G \in \G \otimes I_2} GTG^*  
   = \frac14 T + \frac18 \sum_{G \ne \pm I_4} GTG^*
\end{align*}
It is easy to check that $\Phi(T)$ commutes with $\G \otimes I$,
and hence lies in $\bC I_2 \otimes \B(\H)$.
Hence 
\begin{align*}
  \dist(T, \bC I_2 \otimes \B(\H)) &\le \| T - \Phi(T) \|  
   \le \frac34 \max \| TG - GT \|.
\end{align*}
However 
\begin{align*}
 \| T(2P-1) - (2P-I)T \| &= 2 \| P^\perp TP - PTP^\perp \| \\
 &= 2 \max \{ \| P^\perp TP\|, \|PTP^\perp \| \} .
\end{align*}
Similarly we obtain the same for $G=i(2Q-I)$.
Hence it follows that
\[
 \dist(T, \bC I_2 \otimes \fM_2) \le \frac32 
 \max_{P \in \Lat(\bC I_2 \otimes \fM_2)} \| P^\perp TP\| .
\]

To handle $\bC I \otimes \B(\H)$, form the group
\[
 \G_n = \G^{(n)} \otimes I
  = (\G \otimes \G \otimes \dots \otimes \G) \otimes I
\]
as a subgroup of unitaries in $\fM_{2^n} \otimes \bC I$.
Observe that the commutant is $\bC I_{2^n} \otimes \B(\H)$.
Also each element is the tensor product of elements with
2 point spectrum.  So this property is preserved.
Averaging over $\G_n$ is an expectation $\Phi_n$ onto
$\bC I_{2^n} \otimes \B(\H)$.
As above, we obtain
\begin{align*}
 \dist(T,\bC I_{2^n} \otimes \B(\H)) &\le \| T - \Phi_n(T) \| 
 \le \max_{G \in \G_n} \| TG - GT \| \\
 &\le 2 \max_{P \in \Lat(\bC I_2 \otimes \fM_2)} \| P^\perp TP\| .
\end{align*}
Now one obtains the same estimate for $\bC I \otimes \B(\H)$
by a routine approximation argument.
\qed

If $\fA$ is any unital WOT-closed algebra then $\fA \botimes \B(\H)$ contains two isometries
with orthogonal ranges. Hence we may apply Bercovici's
Theorem~\cite{Ber} to conclude that $\kappa_{\fA\botimes
  \B(\H)}\le3$. 
Here is a more elementary argument that improves on it.

\begin{rem}\label{constant3}
The estimate of a constant 2 cannot be improved by using an
expectation.  
This is because $\|T-\Phi(T)\|$ can be close to 
$2\dist(T,\bC I \otimes \B(\H))$.
Indeed, forgetting about the tensor product with $\B(\H)$,
consider $T = \frac1n \one_n \one_n^*$ where $\one_n$ is the
vector with $n$ ones.  So $T$ is a projection, and its
distance to the scalars is $\frac12$.
However $\Phi(T) = \frac1n I_n$ and 
$\|T-\Phi(T)\| = 1-\frac1n = (2 - \frac2n) \dist(T,\bC I_n)$.
(As well, we know that the off-diagonal projection $\id - \Phi$ has norm
$2 - \frac2n$ on $\fM_n$ \cite{BCD}.)

However the expectation has an advantage.
Let $\fA$ be an arbitrary weak-$*$ closed subspace and 
consider $\bC I \otimes \fA$.
Observe that for each unitary $G \in \B(\H) \otimes I$,
\[
 \dist (GTG^*, \bC I \otimes \fA) = \dist (T, G^*(\bC I \otimes \fA)G)
 = \dist (T, \bC I \otimes \fA) .
\]
So 
\[
 \dist( \Phi(T), \bC I \otimes \fA) \le \dist (T, \bC I \otimes \fA) .
\]
Similarly, if $P \in \Lat(\bC I \otimes \fA)$, then $G^* PG$ is also
in $\Lat(\bC I \otimes \fA)$.  
Thus $\beta_{\bC I \otimes \fA}(GTG^*) = \beta_{\bC I \otimes \fA}(T)$.
Again averaging yields
\[
 \beta_{\bC I \otimes \fA}(\Phi(T)) \le \beta_{\bC I \otimes \fA}(T) .
\]

Now $\Phi(T) = I \otimes T_0$.
It is a well-known argument due to Arveson that if 
$\dist( T_0, \fA) = r$, then there is a weak-$*$ continuous functional
$\psi$ of norm one on $\B(\H)$ annihilating $\fA$ so that 
$\psi(T_0) \approx r$.
The corresponding trace class operator is put into polar
decomposition as $\psi = \sum_{n\ge1} s_n e_n f_n^*$ where $\{e_n\}$
and $\{f_n\}$ are orthonormal and $\sum_{n\ge1} s_n = 1$.
Then $x = \sum_{n\ge1}^\oplus \sqrt{s_n} e_n$
and $y = \sum_{n\ge1}^\oplus \sqrt{s_n} f_n$ are unit vectors in
$l^2 \otimes \H$ such that
\[
 \psi(X) = \ip{ I \otimes X x,y}
\]
for all $X \in \B(\H)$.  The subspace $\ol{\bC I \otimes \fA x}$
is invariant for $\bC I \otimes \fA$ and orthogonal to $y$.
Hence we conclude that
\[
 \beta_{\bC I \otimes \fA}( \Phi(T) ) = \dist( \Phi(T), \bC I \otimes \fA) = r .
\]
Therefore
\[
 \dist(T, \bC I \otimes \fA) \le
 \| T - \Phi(T) \| + \beta_{\bC I \otimes \fA}( \Phi(T) ) \le
 3 \beta_{\bC I \otimes \fA}(T) .
\]
This yields the constant 3 in a more elementary way than by
applying Berco\-vici's Theorem.
\end{rem}

\begin{lem}\label{completeCT}
For $0 \ne T\in \B(\H)$, where $\dim\H \ge3$, the subspace $\bC T$ has
complete  hyperreflexivity constant one if and only if $\rank T = 1$.
\end{lem}

\Prf If $T$ is rank one, then $T =  syx^*$ for unit vectors $x$ and $y$
and non-zero scalar $s$.  
Thus $\bC T$ is a nest bimodule for the nests
$\M = \{ 0, \bC x, \H \}$ and $\N = \{ 0, \bC y, \H \}$.  
Hence $\bC T \otimes \B(\H)$ is also a nest bimodule, and thus has
distance constant one.

If the rank of $T$ is at least two, one cannot put $\bC T \otimes \B(\H)$
into the form of Theorem~\ref{T:1hyp}.
Indeed one can find three orthonormal vectors $x_1,x_2,x_3$ so that
$Tx_i = s_i y_i$ where $y_1,y_2,y_3$ are orthonormal and $s_1s_2 \ne 0$.
The compression of $\bC T \otimes \B(\H)$ to the domain 
$\spn\{x_1,x_2,x_3\} \otimes \H$ and range $\spn\{y_1,y_2,y_3\} \otimes \H$
contains all elements of the form $s_1A \oplus s_2A \oplus s_3A$ for $A
\in \B(\H)$.
If this were 1-hyper\-reflexive, then the form of Theorem~\ref{T:1hyp} would
have a single block.  However the subspace is neither a nest bimodule nor
is it one dimensional.  
So this compression is not 1-hyper\-reflexive.
By Theorem~\ref{compress}, $\bC T \otimes \B(\H)$ is also not
1-hyper\-reflexive.
\qed

\begin{prop}\label{comp1}
A subspace $\S$ of $\B(\H_1,\H_2)$ has complete distance constant one
if and only if there are partitions of the identity of $\H_2$ and $\H_1$
respectively: 
$\C = \{ C_j : j \in \J \}$ and  $\D = \{ D_j : j \in \J \}$,
and for each $j \in \J$, there are
nest bimodules $\fX_j$ of $C_j \B(\H_1,\H_2) D_j$ so that
\[
 \S := \{ T \in \B(\H) : C_jTD_j \in \fX_j \qforal j \in \J \} 
\] 
\end{prop}

\Prf Since $\S$ must be 1-hyper\-reflexive, Theorem~\ref{T:1hyp} yields the
desired form for $\S$ except that the subspaces $\fX_j$ could have atoms
of the bimodules replaced by 1-dimensional subspaces $\bC T_{jk}$.
The compression to this subspace must still be completely
1-hyper\-reflexive.  So by the preceeding lemma, each $T_{jk}$ must be rank
one.  But then, it is easy to see that $\fX_j$ is a nest bimodule.
\qed

The following easy result of Ionascu \cite[Prop.~1.3]{Io} will be useful.

\begin{lem}\label{multiply}
Suppose that $\S$ and $\T$ are subspaces and that $X$ is an invertible
operator such that $\S X = \T$.  Then one subspace is hyperreflexive if and
only if the other is; and the constants are related by 
\[
 \kappa_\T \le \kappa_\S \|X\| \, \| X^{-1} \| \qand
 \kappa_\S \le \kappa_\T \|X\| \, \| X^{-1} \| .
\]
\end{lem}

It is well-known that if $\fD$ is a von Neumann algebra with abelian
commutant, then it is hyperreflexive with constant at most 2.
As well, Rosenoer \cite{Ros} showed with a more sophisticated argument that
this is also true for abelian von Neumann algebras.
However we need a variant which may include a zero summand.
As we have seen, adding a zero summand will generally increase
the distance constant.  For example, $\B(\H) \oplus \B(\H)$ has constant 1
while $\B(\H) \oplus \B(\H) \oplus 0$ does not.
So a modification of the proof is required.

\begin{prop}\label{offdiag}
Let $\{e_i : i \in \I \}$ and $\{f_i : i \in \I \}$ be orthonormal sets
Hilbert spaces $\H_1$ and $\H_2$ respectively.
Let 
\[
 \D = \wot\spn\{f_i e_i^* \otimes \B(\K) : i \in \I \}
\subset \B(\H_1 \otimes \K, \H_2 \otimes \K) .
\]
Then $\D$ is $($completely$)$ hyperreflexive with constant $\kappa_\D \le 2$.

Moreover, suppose that $\S \subset \D$.  
Let $\Phi$ be the contractive expectation of $\B(\H_1 \otimes \K, \H_2
\otimes \K)$ onto $\D$. Then 
\[
 \beta_\S( \Phi(T) ) \le \beta_\S ( T )
 \qforal T \in \B(\H_1 \otimes \K, \H_2 \otimes \K) .
\]
\end{prop}

\Prf For each subset $X$ of $\I$, let $E(X)$ be the projection onto
$\spn\{e_i : i \in X\} \otimes \K$; and similarly let $F(X)$ be the
projection onto $\spn\{f_i : i \in X\} \otimes \K$.
Let $E_i = E(\{i\})$ and $F_i = F(\{i\})$.
Define an expectation $\Phi$ onto $\D$ by putting the standard product
measure $\mu$ on $2^\I$ and integrating:
\[
 \Phi(T) = \int_{2^\I} (F(X) - F(X^c)) T (E(X) - E(X^c)) \,d\mu(X)
\]
for $T \in \B(\H\otimes\K)$.
Clearly this is a completely contractive map.
Observe that $\Phi(T) = F(\I) \Phi(T) E(\I) = \Phi(F(\I) T E(\I) )$.
Moreover for every element $i \in \I$, $F_i \Phi(T) = \Phi(T) E_i$.
From this, it easily follows that 
\[
 \Phi(T) = \sotsum_{i\in \I} F_i \Phi(T) E_i ,
\]
and thus it belongs to $\D$.  
On the other hand, if $D \in \D$, then 
\[ D = (F(X) - F(X^c)) D (E(X) - E(X^c)) \]
for every $X$ and hence $\Phi(D) = D$.
So $\Phi$ is an expectation onto $\D$.

We use the expectation to compute
\begin{align*}
 \dist(T, \D) &\le \| T - \Phi(T) \| \\
 &\le \sup_{X \subset \I} \| T - (F(X) - F(X^c)) T (E(X) - E(X^c)) \| .
\end{align*}
The proof will be completed by bounding this by $2 \beta_\D (T)$.
This is where the proof is a bit trickier than the von Neumann algebra case.

We may write $T$ as a $3\times3$ matrix with respect to the decompositions
of the domain  and range into 
\[
 \H_1 \otimes \K = \ran E(X) \oplus \ran E(X^c) \oplus \ran E(\I)^\perp
\]
and
\[
 \H_2 \otimes \K = \ran F(X) \oplus \ran F(X^c) \oplus \ran F(\I)^\perp 
\]
as $T = \begin{bmatrix} T_{ij} \end{bmatrix}_{i,j = 1}^3$.
Decompose $T - (F(X) - F(X^c)) T (E(X) - E(X^c))$ as
\begin{align*}
 \begin{bmatrix}
 0&2T_{12}&T_{13}\\ 2T_{21}&0&T_{23}\\ T_{31}&T_{32}&T_{33}
\end{bmatrix} 
&= \frac12 \begin{bmatrix}
 0&T_{12}&T_{13}\\ T_{21}&0&0\\ 0&T_{32}&T_{33}
\end{bmatrix}
 + \frac12 \begin{bmatrix}
 0&T_{12}&0\\ T_{21}&0&T_{23}\\ T_{31}&0&T_{33}
\end{bmatrix} \\
&\ + \frac12 \begin{bmatrix}
 0&T_{12}&0\\ T_{21}&0&T_{23}\\ 0&T_{32}&0
\end{bmatrix} 
 + \frac12 \begin{bmatrix}
 0&T_{12}&T_{13}\\T_{21}&0&0\\ T_{31}&0&0
\end{bmatrix}
\end{align*}
Each term is bounded by $\frac12 \beta_\D (T)$.
For example, with 
\[
 P = E(X)^\perp = E(X^c) + E(\I)^\perp ,
\]
$\ran \D P = \ran F(X^c)$; and $\ran \D E(X) = \ran F(X)$.
Hence 
\begin{align*}
 \beta_\D(T) &\ge 
 \max\big\{ \| F(X^c)^\perp T P \|, \, \| F(X)^\perp T E(X) \| \big\} \\
 &\ge 
 \max\big\{ \| F(X^c)^\perp T E(X)^\perp \|, \, \| F(X^c) T E(X) \| \big\} \\
 &= \| F(X^c)^\perp T E(X)^\perp + F(X^c) T E(X) \| \\
 &= \left\| \begin{bmatrix}
 0&T_{12}&T_{13}\\ T_{21}&0&0\\ 0&T_{32}&T_{33} \end{bmatrix} \right\|
\end{align*}
The other three terms are handled similarly.

For the last claim, let $\S$ be any subspace of $\D$.
For each unit vector $x \in \H_1 \otimes \K$, let $Q_x$ be the projection
onto $\ol{\S x}$.
Since 
\[
 \S = (F(X) - F(X^c)) \S (E(X) - E(X^c)) ,
\]
the vector $x' = (E(X) - E(X^c)) x$ has the same range, 
i.e.\ $Q_{x'} = Q_x$; and $\|x'\| \le 1$. So
\begin{align*}
 \beta_\S( \Phi(T) ) &= \sup_{\|x\|=1} \| Q_x^\perp  \Phi(T) x \| \\
 &\le \sup_{\|x\|=1} \sup_{X \subset \I} 
  \| Q_x^\perp  (F(X) - F(X^c)) T (E(X) - E(X^c)) x \| \\
 &= \sup_{\|x\|=1}\sup_{X \subset \I} \|(F(X) - F(X^c)) Q_{x'}^\perp T x' \| \\
 &\le \sup_{\|x'\| \le 1} \| Q_{x'}^\perp T x' \| = \beta_\S(T) .
\end{align*}
\upqed

\begin{thm}\label{1dim_comp}
If $0 \ne T\in \B(\H)$, the subspace $\bC T$ has
complete  hyperreflexivity constant at most $4$.
\end{thm}

\Prf Use the polar decomposition to write $T = UP$ where $P$ is positive.
We may suppose that the spectrum of $P$ is a countable set with $0$ as the
only limit point.
Indeed, suppose that we have established the result for this case.
Without loss of generality, $\|T\|=1$.
Given any $0<r<1$, write $P_r = \sum_{n\ge0} r^n E_P(r^{n+1},r^n]$
where $E_P$ is the spectral measure for $P$; and $T_r = UP_r$.
Then $rP_r \le P \le P_r$ and there is an invertible operator $S_r$ in
$W^*(P)$ so that $P = S_rP_r$ and $rI \le S_r \le I$.
Then 
\[
 \bC T \otimes \B(\K) = \big( \bC T_r \otimes \B(\K) \big)(S_r \otimes I) .
\]
Thus the two spaces have distance constants related by a constant bounded by
$\|S_r\|\,\|S_r^{-1}\| = r^{-1}$.
If each $\bC T_r$ has complete hyperreflexivity constant bounded by 4, then
it follows by letting $r$ tend to $1$ that so does $\bC T$.

Since $P$ has discrete spectrum, it is diagonalizable.  
So select an orthonormal basis $\{ e_n : n \ge 0 \}$ so that 
$Pe_n = \tau_n e_n$ where $\tau_0 = 1 \ge \tau_n$ for all $n\ge1$.
Let $f_n = Ue_n$.
Let $\D$ be the \wot-closed subspace of $\B(\H \otimes \K)$
spanned by $\{ f_n e_n^* \otimes \B(\K) : n \ge 0 \}$.
Then by Proposition~\ref{offdiag}, $\D$ has complete hyperreflexivity 
constant at most 2.
Let $\S = \bC T \otimes \B(\K)$. 
Moreover if $\Phi$ is the expectation onto $\D$ constructed in
Proposition~\ref{offdiag}, we obtain that 
\[
 \beta_\S(\Phi(T)) \le \beta_\S (T) 
 \qforal T \in \B(\H_1 \otimes \K, \H_2 \otimes \K) .
\]

We now consider a relative distance constant for $\S$ within  $\D$.
Note that 
\[
 \S = \bC T \otimes \B(\K) = 
 \Big\{ \sum_{n\ge0} f_n e_n^* \otimes \tau_n A : A \in \B(\K) \Big\} .
\]
An element of the predual $\D_*$ is given by a sequence
$\phi=(\phi_n)$ where $\phi_n$ is a weak-$*$ continuous functional on 
$\bC f_n e_n^* \otimes \B(\K)$ for each $n\ge0$ and 
$\| \phi\| = \sum_{n\ge0} \|\phi_n\| < \infty$.
As usual, we identify each $\phi_n$ with an element of the space $\fS_1$ of
trace class operators on $\K$.
The pre-annihilator $\S_\perp$ intersects $\D_*$ in the set $\A$ of
functionals satisfying $\sum_{n\ge0} \tau_n \phi_n = 0$ considered as 
an absolutely convergent sum in $\fS_1$.

We claim that $\phi \in \A$ may be decomposed as a sum of functionals in $\A$
which have rank at most one in each entry, and have norms summing to at most
$2\|\phi\|$.
Indeed, for each $n \ge 1$, decompose $\phi_n$ using polar decomposition into
a sum of rank one functionals $\rho_{nj}$ so that 
$\|\phi\| = \sum_j \| \rho_{nj} \|$.
Define a functional $\psi_{nj} = (\psi_{nji})_{i\ge0}$ by
\[
 \psi_{njn} = \rho_{nj}, \quad \psi_{nj0} = - \tau_n \rho_{nj} 
 \qand \psi_{nji} = 0 \text{  otherwise.}
\]
Then it is clear that $\psi_{nj} \in \A$.
Moreover
\[
 \sum_{n\ge1}\sum_j \| \psi_{nj} \|
 = \sum_{n\ge1}\sum_j (1 + \tau_n) \| \rho_{ij} \|
 \le 2 \sum_{n\ge1} \| \phi_n \| \le 2 \|\phi\| .
\]
So the sum $\sum_{n\ge1}\sum_j \psi_{nj}$ converges to an element of $\A$.
It is clear that for $n\ge1$, the $n$th component of the sum is just $\phi_n$.
If the zeroth component is $\psi$, then we have
\[
 \psi + \sum_{n\ge1} \tau_n \phi_n = 0 = \phi_0 + \sum_{n\ge1} \tau_n \phi_n .
\]
Hence $\psi = \phi_0$ and this sum is precisely $\phi$.

Now we use the fact that if $\phi = (\phi_n)_{n\ge0}$ in 
$\D_*$ has the property that each $\phi_n$ is rank one, then
there is a rank one functional of the same norm on $\B(\H \otimes \K)$ which
agrees with $\phi$ on $\D$.
Indeed, we may choose vectors $x_n, y_n \in \K$ so that $\phi_n = y_n x_n^*$
and $\|x_n\|_2 = \|y_n\|_2 = \|\phi_n\|^{1/2}$.
Let $x = \sum_{n\ge0} e_n \otimes x_n$ and $y = \sum_{n\ge0} f_n \otimes y_n$.
Then 
\[
 \|x\|^2 = \|y\|^2 = \sum_{n\ge0} \|\phi_n\| = \|\phi\| . 
\]
So $\psi = yx^*$ has $\|\psi\| = \|\phi\|$.
Finally it is evident that the restriction of $\psi$ to $\D$ is equal to
$\phi$.

From the predual formulation of hyperreflexivity, we can conclude that
for all $D \in \D$, 
\[
 \beta_\S(D) = 
 \sup \big\{ |\phi(D)| :
 \phi = (\phi_n) \in \A,\, \|\phi\|=1,\, \rank \phi_n \le 1 \FORAL n \big\}
.\]
The calculation above shows that the convex hull of these rank one 
functionals contains the ball in $\A$ of radius $1/2$.
Hence 
\[
 \dist(D,\S) \le 2 \beta_\S(D) .
\]

Now consider $T \in \B(\H \otimes \K)$.
Then 
\begin{align*}
  \dist(T,\S) &\le \| T - \Phi(T) \| + \dist(\Phi(T),\S) \\
  &\le 2 \beta_\D(T) + 2 \beta_\S(\Phi(T)) \le 4 \beta_\S(T) .
\end{align*}
\upqed

\begin{rem}\label{R:rank2}
We can do better in this analysis if $T$ has rank 2.
In this case, there is no loss of generality in taking
$T_s = \diag(1,s,0,0,\dots)$ where $0 < s \le 1$.
Then for an infinite dimension Hilbert space $\K$,
the space $\S = \bC T_s \otimes \B(\K)$ is unitarily equivalent
to the subset of $\fM_3(\B(\K))$ given by 
\[
 \left\{ 
 \begin{bmatrix}A&0&0\\0&sA&0\\0&0&0\end{bmatrix} : A \in \B(\K) 
 \right\} .
\]

There are now two improvements in the argument of Theorem~\ref{1dim_comp}.
The first is that the expectation $\Phi$ onto $\fD_3 \otimes \B(\K)$
yields a better estimate because one averages over the finite
group of diagonal matrices with $\pm1$ as entries.  
As in the proof of Proposition~\ref{P:scalartensor}, there is an economy
because two of the 8 group elements are $\pm I$.  
So one obtains an upper bound of 
$\frac32 \beta_{\fD_3 \otimes \B(\K)}(X)$ for
$\dist(X, \fD_3 \otimes \B(\K))$.  
Actually, the distance constant for $\fD_3 \otimes \B(\K)$ is known to be
exactly $\sqrt{3/2}$ \cite{DO}.  
However the expectation has the advantage that
$\beta_\S( \Phi(X) ) \le \beta_\S(X)$, which we do not know for the
closest point.

The second improvement is that $\S$ has a relative distance constant of 1
within $\fD_3 \otimes \B(\K)$.
This can also be seen from the proof of the previous theorem.
Indeed, $\S_\perp \cap \big( \fD_3 \otimes \B(\K) \big)_*$ consists
of $\phi = (\phi_n)_{n\ge0}$ such that $\phi_{00} = -s \phi_{11}$.
Decompose each $\phi_n$ for $n \ge 1$ into a sum of rank one elements
$\rho_{nj}$ so that $\| \phi_n\| = \sum_j \|\rho_{nj} \|$.
Then set $\psi_{1j} = (-s \rho_{1j}, \rho_{1j},0,0,\dots)$;
and
set $\psi_{nj}$ to have $\rho_{nj}$ in the $n$th entry and 0 elsewhere.
Then each $\psi_{nj}$ belongs to 
$\S_\perp \cap \big( \fD_3 \otimes \B(\K) \big)_*$ and is rank at most
one in each entry.
Moreover the norms sum exactly to $\|\phi\|$.
The proof is completed as above.
So one obtains a distance constant of at most $1.5+1 = 2.5$.
\end{rem}

The distance constant fails to be continuous except in rare cases
\cite{PV,Io}.  The example in the remark above displays this
in a striking way that we have not seen before.
Observe that Proposition~\ref{multiply} shows that the constant
$\kappa_{\bC T_s}^c$ is a continuous function of $s$ for $s > 0$
\cite{Io}.

\begin{prop}\label{1dim_small}
Let $T_s = \diag(1,s,0,0,0 \dots)$ for $0 \le s \le 1$.
Then $\sqrt 2 \le \dlim_{s \to 0^+} \kappa_{\bC T_s}^c \le \frac52$
while $\kappa_{\bC T_0}^c = 1$.
\end{prop}

\Prf  The upper bound of 2.5 was just established.
And $\kappa_{\bC T_0}^c = 1$ follows from Lemma~\ref{completeCT}.
As above, we consider $\S_s = \bC T_s \otimes \B(\K)$ as
the set of operators of the form
$\begin{bmatrix}X&0&0\\0&sX&0\\0&0&0\end{bmatrix}$ for $X \in \B(\K)$.

We observe as in Theorem~\ref{1dim_comp} that 
$( \S_s )_\perp$ consists of all trace class 
operators  $\phi = \begin{bmatrix} \phi_{ij}\end{bmatrix}_{i,j=1}^3$
such that $\phi_{11} = -s \phi_{22}$.
Let $e_1$ and $e_2$ be two orthonormal vectors and $f$ any unit vector,
and consider
\[
 \psi = \begin{bmatrix}
         \frac{-s}{s+\sqrt2} fe_2^* & 0 & 0\\
         \frac{1}{s+\sqrt2} fe_1^* & \frac{1}{s+\sqrt2} fe_2^*&0\\
         0&0&0\end{bmatrix} 
 \qand
 A = \begin{bmatrix} -e_2 f^* & \frac1{\sqrt2} e_1 f^* &0\\
       0 & \frac1{\sqrt2} e_2 f^* &0 \\ 0&0&0 \end{bmatrix}.
\]
It is easy to verify that $\psi \in ( \S_s )_\perp$ and 
$\| \psi \|_1 = 1$; and that $\|A\| = 1 = \psi(A)$. 
So $\dist(A, \S_s) = \|A\|$.  
Hence $\kappa_{\S_s} \ge \beta_{\S_s}(A)^{-1}$.

The predual formulation of the constant shows that $\beta_{\S_s}(A)$
is obtained as $\sup | \phi(A) |$ taken over all rank one elements of the
unit ball of $( \S_s )_\perp$.
Note that the compression of any such functional to the upper left
$2\times2$ corner still lies in $( \S_s )_\perp$ and still has rank one.
As $A$ is supported in this $2\times2$ corner, the set of
compressions will yield the same supremum.
This reduces the problem to $2\times2$ matrices.
That is,
\[
 A = \begin{bmatrix} -e_2 f^* & \frac1{\sqrt2} e_1 f^*\\
       0 & \frac1{\sqrt2} e_2 f^* \end{bmatrix}
 \AND
 (\S_s)_\perp = \left\{ \begin{bmatrix}
 -s \phi_{22} & \phi_{12} \\ \phi_{21} & \phi_{22}\end{bmatrix} : 
 \phi_{ij} \in \fS_1(\K) \right\} .
\]

Here the description of rank one elements of $(\S_s)_\perp$ is
particularly easy. There are three families:
\[
 \begin{bmatrix}0&0\\xy^*&0\end{bmatrix}\,, \qquad
 \begin{bmatrix}0&xy^*\\0&0\end{bmatrix} \qand
 \begin{bmatrix}-saxy^* & b xy^* \\ cxy^* & a xy^* \end{bmatrix} 
\]
where $\|x\| = \|y\| \le 1$ and $a,b,c \in \bC$ satisfy
\[
 sa^2 + bc = 0 \qand (1+s^2)|a|^2 + |b|^2 + |c|^2 \le 1 .
\]
Indeed, a rank one that has 0 in the $2,2$ entry must be supported in
either the $2,1$ entry or the $1,2$ entry.
When $\phi_{22} = axy^*$ for $\|x\|=\|y\|=1$, this forces
$\phi_{11} = -saxy^*$.  Then to make $\phi$ rank one, the other two
entries must also be multiples of $xy^*$.
The condition $sa^2 + bc = 0$ is a determinant condition equivalent to
being rank one.  
Then the trace norm equals the Hilbert-Schmidt norm, so it is easily
calculated and one obtains $(1+s^2)|a|^2 + |b|^2 + |c|^2 \le 1$.

Evaluating the first two classes on $A$ yields $1/\sqrt2$ and $0$
respectively.  So consider the third class.  
Then from Cauchy-Schwarz,
\begin{align*}
 \sup |\phi(A)| &=
 \sup \big| (s + \tfrac1{\sqrt2}) a \ip{e_2,y} \ip{x,f}
   + \frac a {\sqrt2} \ip{e_2,y} \ip{x,f} \big| \\
 &= \sup \big( (s + \tfrac1{\sqrt2})^2 |a|^2 + \tfrac12 |c|^2 \big)^{1/2}
 \\ &= \sup \sqrt{ \frac{(2s^2+2\sqrt2 s + 1) |a|^2 + |c|^2 } 2 }
\end{align*}
The constraints yield $b = -sa^2/c$ and
\[
 1 \ge (1+s^2)|a|^2 + |b|^2 + |c|^2 = 
 \frac{ (1+s^2)|a|^2|c|^2 + s^2 |a|^4 + |c|^4 }{|c|^2} .
\]

This problem may be solved by Lagrange multipliers.  
Since we are interested in the limit of this supremum as $s$ tends to
$0^+$, we may observe that this will be the solution to the extremal
problem
\[
 \sup \Big\{ \sqrt{\frac{x^2+y^2}2} : 0 < (x^2 + y^2)y^2 \le y^2 \Big\} 
 = \frac1{\sqrt2} .
\]
\upqed


\end{document}